\documentclass[a4paper,reqno,oneside]{amsart}
\usepackage[all,cmtip]{xy}
\newdir{ >}{{}*!/-5pt/@{>}}
\usepackage{amssymb,latexsym,amsmath}
\usepackage{amsthm}
\usepackage[a4paper]{geometry}
\usepackage[colorlinks=true,breaklinks=true,linkcolor=black,citecolor=black,hidelinks]{hyperref}
\usepackage{verbatim}
\usepackage[textsize=footnotesize]{todonotes}
\usepackage{times}
\usepackage{manfnt}
\usepackage{mathrsfs}
\theoremstyle{plain}

\swapnumbers
\theoremstyle{plain}
\newtheorem{theorem}{Theorem}[section]

\newtheorem{lemma}[theorem]{Lemma}
\newtheorem{prop}[theorem]{Proposition}
\newtheorem{corollary}[theorem]{Corollary}

\theoremstyle{definition}
\newtheorem{ntheorem}[theorem]{\theoremname}
\newcommand{\theoremname}{Test}
\newenvironment{des}[1]{\renewcommand{\theoremname}{#1} \begin{ntheorem}}{\end{ntheorem}}
\newcommand{\bdes}{\begin{des}}
\newcommand{\edes}{\end{des}}
\newtheorem{ex}[theorem]{Example}
\newtheorem{definition}[theorem]{Definition}
\newtheorem{exer}{Exercise}
\newtheorem{rem}[theorem]{Remark}
\newtheorem{observation}[theorem]{Observation}
\newtheorem{app}[theorem]{Application}
\newtheorem{his}[theorem]{Historical Note}
\newtheorem*{notation}{Notation}
\newtheorem*{conv}{Convention}
\newtheorem{Observation}[theorem]{Observation}
\newtheorem{cons}[theorem]{Construction}
\newtheorem{diss}[theorem]{}
\renewcommand{\theenumi}{\roman{enumi}}

\newcommand{\bthm}{\begin{theorem}}
\newcommand{\ethm}{\end{theorem}}

\newcommand{\bprop}{\begin{prop}}
\newcommand{\eprop}{\end{prop}}

\newcommand{\bcoro}{\begin{corollary}}
\newcommand{\ecoro}{\end{corollary}}

\newcommand{\bex}{\begin{ex}}
\newcommand{\eex}{ \end{ex}}

\newcommand{\bdf}{\begin{definition}}
\newcommand{\edf}{ \end{definition}}

\newcommand{\brem}{\begin{rem}}
\newcommand{\erem}{ \end{rem}}

\newcommand{\blem}{\begin{lemma}}
\newcommand{\elem}{\end{lemma}}

\newcommand{\bhis}{\begin{his}}
\newcommand{\ehis}{ \end{his}}

\newcommand{\bobs}{\begin{observation}}
\newcommand{\eobs}{ \end{observation}}

\newcommand{\bproof}{\begin{proof}}
\newcommand{\eproof}{\end{proof}}

\newcommand{\bexr}{\begin{exer}}
\newcommand{\eexr}{ \end{exer}}

\newcommand{\bnotation}{\begin{notation}}
\newcommand{\enotation}{\end{notation}}

\newcommand{\bconv}{\begin{conv}}
\newcommand{\econv}{ \end{conv}}

\newcommand{\bapp}{\begin{app}}
\newcommand{\eapp}{ \end{app}}

\newcommand{\bfact}{\begin{Observation}}
\newcommand{\efact}{ \end{Observation}}

\newcommand{\bcons}{\begin{cons}}
\newcommand{\econs}{ \end{cons}}

\newcommand{\bdiss}{\begin{diss}}
\newcommand{\ediss}{ \end{diss}}

\newcommand{\Mod}{\mathsf{Mod}}
\renewcommand{\mod}{\mathsf{mod}}

\newcommand{\Add}{\mathsf{Add}}
\newcommand{\add}{\mathsf{add}}

\newcommand{\Gen}{\mathsf{Gen}}

\newcommand{\Proj}{\mathcal{P}}

\newcommand{\GP}{\mathcal{GP}}
\newcommand{\Gp}{\mathcal{GP}^\mathrm{fin}}

\newcommand{\fin}{\mathrm{fin}}

\newcommand{\arrow}{\longrightarrow}
\renewcommand{\mapsto}{\longmapsto}

\renewcommand{\tilde}{\widetilde}
\renewcommand{\hom}{\operatorname{Hom}}

\newcommand{\Ext}{\operatorname{Ext}}
\newcommand{\Gext}{\operatorname{Gext}}
\newcommand{\ext}{\operatorname{ext}}

\newcommand{\dirlim}{\varinjlim  \limits}

\renewcommand{\bar}{\overline}
\newcommand{\End}{\operatorname{End}}

\renewcommand{\ker}{\operatorname{Ker}}
\newcommand{\coker}{\operatorname{Coker}}
\newcommand{\im}{\operatorname{Im}}

\newcommand{\FPD}{\operatorname{FPD}}
\newcommand{\fpd}{\operatorname{fpd}}
\newcommand{\pd}{\operatorname{pd}}

\newcommand{\id}{\operatorname{id}}

\newcommand{\Gpd}{\operatorname{Gpd}}

\newcommand{\op}{{\operatorname{op}}}

\newcommand{\iso}{\mathop{\cong} \limits}
\newcommand{\hookarrow}{\xymatrix{  \ar@{^{(}->}[r] & }}
\newcommand{\bperp}{{\perp\!\!\frak B}}
\newcommand{\Gperp}{\mathrm{G}\!\perp}
\newcommand{\xpd}{\operatorname{\mathcal{X}-pd}}
\newcommand{\yid}{\operatorname{\mathcal{Y}-id}}
\newcommand{\Ab}{\mathsf{Ab}}
\newcommand{\GI}{\mathcal{GI}}

\title{\textbf{Infinitely Generated Gorenstein Tilting Modules}}
\author{Pooyan Moradifar \and Siamak Yassemi}
\address[Pooyan Moradifar]{School of Mathematics, Statistics and Computer Science \\ University of Tehran \\ Tehran \\ Iran}
\email{pmoradifar@ut.ac.ir}
\address[Siamak Yassemi]{School of Mathematics, Statistics and Computer Science \\ University of Tehran \\ Tehran \\ Iran}
\email{yassemi@ut.ac.ir}
\date{\today}
\keywords{$\hom$-balanced pair; Cotorsion pair; Gorenstein tilting module; Gorenstein tilting class; Finitistic dimensions}

\begin{document}

\begin{abstract}
The theory of finitely generated relative (co)tilting modules has been established in the 1980s by Auslander and Solberg, and 
infinitely generated relative tilting modules have recently been studied by many authors in the context of Gorenstein homological algebra.
In this work, we build on the theory of infinitely generated Gorenstein tilting modules by developing ``Gorenstein tilting approximations''
and employing these approximations to study Gorenstein tilting classes and their associated relative cotorsion pairs. As applications of our results, 
we discuss the problem of existence of complements to partial Gorenstein tilting modules as well as some connections between Gorenstein tilting modules and 
finitistic dimension conjectures.
\end{abstract}
\maketitle


\section*{Introduction}

Tilting theory has been developed originally in the 1980s through the work of Brenner, Butler~\cite{brenner.butler}, Happel, Ringel~\cite{tiltedalg.happel,happel.book}, and Miyashita~\cite{miyashita},  and it soon became a thriving branch of representation theory of algebras with invaluable applications; cf.~\cite{tilting.handbook}.
Since its advent, the theory has been generalized in various directions and contexts. The two generalizations which are of particular interest in the present work
are:
\begin{itemize}
\item Generalization of tilting theory to the context of infinitely generated modules over arbitrary rings by Colpi, Trlifaj~\cite{inf.tilting}, Angeleri-H\"{u}gel and Coelho~\cite{hugel.coel};

\item Generalization of tilting theory to relative setting by Auslander and Solberg~\cite{ASII}.
\end{itemize}
The scope of the present work lies within the intersection of  the above-mentioned branches of tilting theory, namely
``infinitely generated relative tilting theory''.

Relative homological algebra originates in the work of Eilenberg and Moore~\cite{EM} in the 1960s, and has been revived through the theory of approximation of modules developed by Auslander et al.~\cite{AR.applications,aus.preproj} and Enochs et al.~\cite{enochs.coversenv,enochs.book}.
In the 1990s, Auslander and Solberg~\cite{ASI,ASII,ASIII} initiated a systematic study of relative homological algebra  in representation theory of algebras,
and in particular they introduced in~\cite{ASII} the notion of a ``relative (co)tilting module'' and
generalized many results of the standard (co)tilting theory to the relative setting. However, the scope of the relative (co)tilting theory developed by Auslander and
Solberg was limited to finitely generated modules over artin algebras, and therefore attempts have been made
in recent years  to transcend the scope of the theory to the context of infinitely generated modules, parallel to the theory of (standard) infinitely generated tilting modules. This line of thought has been followed especially in the context of ``Gorenstein homological algebra''---a particular relative (co)homology theory---by several authors, e.g.~\cite{wei},~\cite{Gtilting}, and~\cite{tilting.sub}, where the theory of ``infinitely generated Gorenstein tilting modules''
has been developed to some extent.

Tools of approximation theory of modules play a significant role in studying infinitely generated tilting modules; cf.~\cite{trlifaj.infinite.cotorsion,GT}.
One of the missing tools in studying infinitely generated Gorenstein tilting modules---in comparison to the (standard) infinitely generated tilting modules---is  ``Gorenstein tilting approximations''. The main goal of the present paper is to develop/sharpen the tools of approximation theory for studying
infinitely generated Goresntein tilting modules. Thus, the structure of the paper is as follows: 

In Section~1, we recall some notions of relative homological algebra which will be used in the subsequent sections. The key notions of this section are 
`` $\hom$-balanced pairs'' and their associated ``relative Ext-functors''.

In Section~2, we study ``relative orthogonal classes'', which are classes of modules defined as roots of relative Ext-functors. The key notion of this section is the notion of ``proper filtration''~\ref{df:rel.filtration} and the key result is the relative Eklof Lemma~\ref{Eklof.lemma.Rel}, which is the relative analogue of
the classical Eklof Lemma; cf.~\cite[1.2]{eklof} and~\cite[Lemma~6.2]{GT}.

Section~3 contains fundamental results which provide us with ``sharp'' approximations suitable for studying infinitely Gorenstein tilting modules in the 
subsequent sections. The key notion of this section is the notion of a ``relative cotorsion pair'' associated with a $\hom$-balanced pair, and the 
the main result of this section is Theorem~\ref{ET.rem} about completeness of these relative cotorsion pairs, which is the relative counterpart of the Eklof-Trlifaj's completeness theorem for (standard) cotorsion pairs; cf.~\cite[Theorem~10]{trlifaj.Extvanish} and~\cite[Theorem 6.11]{GT}.

With the sharp approximation tools developed in Section~3 at our disposal, we embark on investigating infinitely generated Gorenstein tilting modules
in Section~4  and Section~5. General results about infinitely generated Gorenstein tilting modules and their approximations are discussed in Section~4.
These results pave the way for the results of Section~5, where we focus on infinitely generated Gorenstein tilting modules over virtually Gorenstein CM-finite rings.
The main result of this section is Theorem~\ref{thm:Gtilting.class} which characterizes Gorenstein tilting classes, parallel to Angeleri-H\"{u}gel and Coelho's
characterization~\cite[Theorem~4.1]{hugel.coel} of tilting classes, and allows us to tie the theory of infinitely generated Gorenstein tilting modules to finitistic dimension conjectures in the next section.

Some applications of infinitely generated Gorenstein tilting modules were discuss in Section~6. In the first part of this section we discuss the problem of existence of
complements to partial Gorenstein tilting modules. It is proved in Theorem~\ref{thm:partial.Gtilting} that every partial Gorenstein tilting module can be ``completed'' to a Gorenstein tilting module. This result is the Gorenstein analogue of a result of Angeleri-H\"{u}gel and Coelho~\cite[Theorem~2.1]{hugel.coel.partial} and can be regarded as a non-finitely generated version of~\cite[Proposition~3.25]{ASII} in the Gorenstein setting.
In the rest of the section, we study some connections between infinitely generated Gorenstein tilting modules and finitistic dimension conjectures,
parallel to~\cite{trlifaj.fpd}.

\section{Preliminaries}
\label{sec:pre}

In this preliminary section, we recall some basic notions and facts from relative homological algebra which will be used in the subsequent sections.
For unexplained notions of relative homological algebra we refer the reader to~\cite{enochs.book} and~\cite{GT}.

\bconv
Throughout this paper,  a ``ring'' means an arbitrary ring with non-zero unit.
Such a ring will be denoted by $\Lambda$. All modules are \textit{left modules}, and right modules over $\Lambda$ will be considered as left modules 
over $\Lambda^\op$---the opposite ring of $\Lambda$.  
\econv

\bnotation
The category of $\Lambda$-modules is denoted by $\Mod (\Lambda)$ and the category of 
\emph{strongly finitely presented} $\Lambda$-modules, i.e. $\Lambda$-modules admitting a degreewise finitely generated projective resolution, is denoted by $\mod (\Lambda)$. We also let $\Ab := \Mod (\mathbb{Z})$.

For any class $\mathcal{C}$ of $\Lambda$-modules we let $\mathcal{C}^\fin :=\mathcal{C} \cap \mod (\Lambda)$ and we denote by $\dirlim \mathcal{C}$ the class of all $\Lambda$-modules which are the direct limit of a direct system of $\Lambda$-modules in $\mathcal{C}$.

For any integer $n \geq 0$, the class of $\Lambda$-modules of projective dimension at most $n$ will be denoted by $\Proj_n$ and
the class of $\Lambda$-modules of injective dimension at most $n$ will be denoted by $\mathcal{I}_n$. 
\enotation

\bdes{Approximations} \label{approximations}
Let $\mathcal{X}$ be a class of $\Lambda$-modules. A $\Lambda$-complex
\[\xymatrix{C_\bullet =  \cdots \ar[r] & C_{n+1} \ar[r] & C_{n} \ar[r] & C_{n-1} \ar[r] & \cdots }\]
is said to be \emph{$\hom_\Lambda (\mathcal{X} , -)$-exact} (respectively, \emph{$\hom_\Lambda (-,\mathcal{X})$-exact})
if for every $X \in \mathcal{X}$ the induced complex $\hom_\Lambda (X, C_\bullet)$ (respectively, $\hom_\Lambda (C_\bullet , X)$)
is exact.

Recall that a $\Lambda$-homomorphism $f : X \arrow M$ with $X \in \mathcal{X}$ is called an \emph{$\mathcal{X}$-precover} 
of $M$ in $\Mod (\Lambda)$ if the $\Lambda$-complex $\xymatrix{X \ar[r]^-f & M \ar[r] & o}$ is $\hom_\Lambda (\mathcal{X} , -)$-exact. The class $\mathcal{X}$ 
is called a \emph{precovering class} in $\Mod (\Lambda)$ provided that each $\Lambda$-module admits an $\mathcal{X}$-precover in $\Mod (\Lambda)$.
The notions of a \emph{preenvelope} and a \emph{preenveloping class} are defined dually; see e.g.~\cite[Section~7]{enochs.book}. 
If we restrict the definition of a ``precovering class'' and a ``preenveloping class'' to the realm of $\mod (\Lambda)$, then
we recover the notion  of a \emph{contravariantly finite} and \emph{covariantly finie} class of modules respectively, as in~\cite{aus.preproj} and~\cite{AR.applications}.
\edes

\bdes{Relative Resolutions and Syzygies} \label{rel.resolution} \label{rel.syzygy}
Let $\mathcal{X}$ be a class of $\Lambda$-modules.
Given a $\Lambda$-module $M$, a $\hom_\Lambda (\mathcal{X} , -)$-exact $\Lambda$-complex of the form
\[\xymatrix{X_\bullet = \cdots \ar[r] & X_1 \ar[r]^-{d_1} & X_0 \ar[r]^-{d_0} & M \ar[r] & o}\]
wherein each $X_i$ belongs to $\mathcal{X}$ is called a \emph{left $\mathcal{X}$-resolution} of $M$. If we delete the map $d_0$ from $X_\bullet$, 
then the resulting complex is called a \textit{deleted left $\mathcal{X}$-resolution} of $M$.
Furthermore, for every integer $i \geq 0$ the $\Lambda$-module $\im (d_i)$ is called the \textit{$(i-1)$-th syzygy module} of $M$ relative to $X_\bullet$. The notions of a \emph{(deleted) right $\mathcal{X}$-resolution} and its respective \emph{cosyzygies} are defined dually; see e.g.~\cite[8.1]{enochs.book}. It is easily seen (see~\cite[Proposition 8.1.3]{enochs.book}) that every module admits a left $\mathcal{X}$-resolution
(respectively, right $\mathcal{X}$-resolution) if and only if the class $\mathcal{X}$ is precovering (respectively, preenveloping).
\edes

\bdes{Relative Dimensions} \label{rel.dim}
Let $\mathcal{X}$ and $\mathcal{Y}$ be full additive subcategories of $\Mod (\Lambda)$ where $\mathcal{X}$ is 
precovering and $\mathcal{Y}$ is preenveloping. The \emph{$\mathcal{X}$-projective dimension} \index{$\mathcal{X}$-projective dimension} of
a $\Lambda$-module $M$, denoted by $\xpd_R (M)$, is defined as the minimum length of  left $\mathcal{X}$-resolutions of $M$. That is,
\[\xpd_R (M ) = \inf \Big \{ \sup \left \{ n \geq 0 : X_n \ne o \right \}  \mid 
\text{$X_\bullet$ is a left $\mathcal{X}$-resolution of $M$} \Big \} \: . \]
Note that for $\mathcal{X}=\mathcal{P}_0$ we have $\xpd_\Lambda (-)=\pd_\Lambda (-)$, the usual projective dimension in homological algebra. Dually, the \emph{$\mathcal{Y}$-injective dimension} of $M$, denoted by $\yid_\Lambda (M)$, is defined as the minimum length of right $\mathcal{Y}$-resolutions of $M$. That is,
\[\yid_R (M ) = \inf \Big \{ \sup \left \{ n \geq 0 : X_{-n} \ne o \right \} \mid
\text{$Y_\bullet$ is a right $\mathcal{Y}$-resolution of $M$} \Big \} \: . \]
Note that for $\mathcal{Y}=\mathcal{I}_0$ we have $\yid_\Lambda (-)=\id_\Lambda (-)$, the usual injective dimension in homological algebra.
\edes

\bdes{Relative Ext-functors} \label{rel.Ext}
Let $(\mathcal{X} , \mathcal{Y})$ be a pair of full additive subcategories of $\Mod (\Lambda)$ where $\mathcal{X}$ is precovering and $\mathcal{Y}$ is preenveloping.
For $\Lambda$-modules $M$ and $N$ and every integer $n$, let
\[ \Ext^n_\mathcal{X} (M,-)  : \Mod (\Lambda) \arrow \Ab \] denote the \textit{$n$-th relative right derived functor} of the covariant 
functor $\hom_\Lambda (M,-)$, and also let 
\[\ext^n_\mathcal{Y} (-,N)  : \Mod (\Lambda)^\op \arrow \Ab\]
denote the \textit{$n$-th relative right derived functor} of the contravariant functor $\hom_\Lambda (-,N)$; see~\cite[Section~8.2]{enochs.book} for more details. Thus, if $X_M$ is a deleted left $\mathcal{X}$-resolution of $M$ and $Y^N$ is a deleted right $\mathcal{Y}$-resolution of $N$, then
\begin{align*}
\Ext^n_\mathcal{X} (M,N) & = H_{-n} \big ( \hom_\Lambda (X_M , N) \big ) \: , \\
\ext^n_\mathcal{Y} (M,N) & =  H_{-n} \big ( \hom_\Lambda (M , Y^N) \big ) \: ,
\end{align*}
for any $n \in \mathbb{Z}$. In the special case where $\mathcal{X}=\mathcal{P}_0$ and $\mathcal{Y}=\mathcal{I}_0$, it is well-known (see e.g.~\cite[Theorem 6.67]{RotmanSE})
that the bi-functors
\begin{align*}
\Ext^n_{\mathcal{P}_0} (-,-) : \Mod (\Lambda)^\op \times \Mod (\Lambda) \arrow \Ab \: , \\
\ext^n_{\mathcal{I}_0} (-,-) : \Mod (\Lambda)^\op \times \Mod (\Lambda) \arrow \Ab \: ,
\end{align*}
are isomorphic, and they are denoted commonly  by $\Ext^n_\Lambda (-,-)$, which is the usual Ext-functor in homological algebra. This property holds essentially because $\big (\mathcal{P}_0 , \mathcal{I}_0 \big )$ is a ``$\hom$-balanced pair'' in the sense defined below.
\edes

\bdes{$\mathbf{Hom}$-balanced Pairs} \label{df:balanced.pair}
Let $\mathfrak{B} : = (\mathcal{X} , \mathcal{Y})$ be a pair of full additive subcategories of $\Mod (\Lambda)$ which are closed under direct summand, and
assume that $\mathcal{X}$ is precovering and $\mathcal{Y}$ is preenveloping. The pair $\mathfrak{B}$ is  
called a \emph{$\hom$-balanced pair} in $\Mod (\Lambda)$ if $\hom_\Lambda (-,-)$ is 
``right balanced'' by $\mathfrak{B}$ in the sense of~\cite[Definition 8.2.13]{enochs.book}.
That is, for every $\Lambda$-module $M$, there is a left $\mathcal{X}$-resolution of $M$ which is $\hom_\Lambda (-,\mathcal{Y})$-exact and there is 
a right $\mathcal{Y}$-resolution of $M$ which is $\hom_\Lambda (\mathcal{X} ,  -)$-exact. $\hom$-balanced pairs have been studied (sometimes under 
different names) by many authors; see e.g.~\cite{enochs.balanced},~\cite{enochs.Gext},~\cite{chen.balanced} and~\cite{enochs.Fbalance}.
\edes

\bdes{$\mathfrak{B}$-exactness} \label{df:B-exactness}
Given a $\hom$-balanced pair $\frak B : = (\mathcal{X} , \mathcal{Y})$ in $\Mod (\Lambda)$, it is straightforward to see (cf.~\cite[Proposition~2.2]{chen.balanced}) 
that a $\Lambda$-complex $C_\bullet$ is $\hom_\Lambda (\mathcal{X} ,- )$-exact if and only if it is $\hom_\Lambda (-,\mathcal{Y})$-exact.
In this case, the $\Lambda$-complex $C_\bullet$ is said to be  \emph{$\frak B$-exact}.
A $\frak B$-exact complex of $\Lambda$-modules of the form
\[\xymatrix{ o \ar[r] & M' \ar[r]^-f & M \ar[r]^-g & M'' \ar[r] & o }\]
is called a \emph{short $\frak B$-exact sequence} or a \emph{$\frak B$-extension} of $M'$ by $M''$.
In this case, $f$ is called a \textit{$\mathfrak{B}$-monomorphism} and $g$ is called a\textit{ $\mathfrak{B}$-epimorphism}.
\edes

Note that $\frak B$-exact complexes are \emph{not} necessarily exact in general,
but they will be exact if the $\hom$-balanced pair $\frak B$ is ``admissible'' in the sense defined below.

\bdes{Admissible $\mathbf{Hom}$-balanced Pairs} \label{df:admissible.Bpair}
A  $\hom$-balanced pair $\frak B := (\mathcal{X} , \mathcal{Y})$ in $\Mod (\Lambda)$ is said to
be \emph{admissible} if one of the following equivalent conditions holds (cf.~\cite[Proposition~2.2 and Corollary~2.3]{chen.balanced}):
\begin{enumerate}
\item $\mathcal{X}$-precovers are surjective;

\item Every $\hom_\Lambda (\mathcal{X},  -)$-exact complex of $\Lambda$-modules is exact;

\item Every $\hom_\Lambda (-,\mathcal{Y})$-exact complex of $\Lambda$-modules is exact;

\item $\mathcal{Y}$-preenvelopes are injective.
\end{enumerate}
Thus, if $\frak B$ is admissible, then $\frak B$-exact complexes are exact in the usual sense. 
\edes

\bobs \label{pullback.of.Bexact}
Let $\frak B$ be an admissible $\hom$-balanced pair in $\Mod (\Lambda)$. Given a $\frak B$-exact sequence
\[\xymatrix{ \delta := o \ar[r] & M' \ar[r]^-f & M \ar[r]^-g & M'' \ar[r] & o }\]
of $\Lambda$-modules and $\Lambda$-homomorphisms $\alpha : N'' \arrow M''$ and $\beta : M' \arrow N'$, we may form the pullback diagram
\[\xymatrix{
\delta \cdot \alpha =  o \ar@{..>}[r] & M' \ar@{:}[d]_-{1_{M'}} \ar@{..>}[r] & P \ar@{..>}[r] \ar@{..>}[d] & N'' \ar@{..>}[r] \ar[d]^-\alpha & o \\
\quad \delta =  o \ar[r] & M' \ar[r]^-f & M \ar[r]^-g & M'' \ar[r] & o
}\]
and the pushout diagram
\[\begin{array}{c}
\xymatrix{
\quad \delta =  o \ar[r] & M' \ar[r]^-f \ar[d]_-\beta & M \ar[r]^-g \ar@{..>}[d] & M'' \ar[r] \ar@{:}[d]^-{1_{M''}} & o \\
\beta \cdot \delta =  o \ar@{..>}[r] & N' \ar@{..>}[r] & D \ar@{..>}[r] & M'' \ar@{..>}[r] & o
}\end{array} \: .
\]
It then follows easily from the universal properties of pullback and pushout (see also~\cite[Lemma~2.3]{rel.cotorsion}) that the exact
sequences $\delta \cdot \alpha$ and $\beta \cdot \delta$ are $\frak B$-exact. This observation will be used frequently to construct a new
short $\frak B$-exact sequence from a given short $\frak B$-exact sequence.
\eobs

\bdes{Relative $\mathbf{Ext}$-functors} \label{BExt}
Given a $\hom$-balanced pair $\mathfrak{B} : = (\mathcal{X} , \mathcal{Y})$ in $\Mod (\Lambda)$, one can prove  that
the bi-functors $\Ext^n_\mathcal{X} (-,-)$ and $\ext^n_\mathcal{Y} (-,-)$ are isomorphic; see e.g.~\cite[Theorem 8.2.14]{enochs.book}. In this case, 
the two functors will commonly be denoted by $\Ext^n_{\frak B} (-,-)$. This functor is called the \emph{relative Ext-functor associated with (the $\hom$-balanced pair) $\frak B$}. Thus, given $\Lambda$-modules $M$ and $N$, we have
\[ H_{-n} \big ( \hom_\Lambda (X_M , N) \big ) = \Ext^n_{\frak B} (M,N) = H_{-n} \big ( \hom_\Lambda (M,Y^N) \big ) \: , \]
for every integer $n$, where $X_M$ is a deleted left $\mathcal{X}$-resolution of $M$ and $Y^N$ is a deleted right $\mathcal{Y}$-resolution of $N$.
\edes

\bdes{The Standard $\mathbf{Hom}$-Balanced Pair} \label{standard.balanced.pair}
It is clear from the definition that $\frak S : =  ( \Proj , \mathcal{I}  )$ is an admissible $\hom$-balanced pair, 
which is called the \emph{standard $\hom$-balanced pair}. In this case, $\mathfrak{S}$-exact complexes coincide with exact sequences and
$\Ext^n_{\frak S} (-,-)$ is the standard $\Ext$-functor, $\Ext^n_\Lambda (-,-)$, in homological algebra.
\edes

\bobs \label{obs:BExt} 
If $\frak B : = (\mathcal{X} , \mathcal{Y})$ is a $\hom$-balanced pair in $\Mod (\Lambda)$, then the relative $\Ext$-functor $\Ext^n_{\frak B} (-,-)$ 
possesses all the properties of the functors $\Ext^n_\mathcal{X}$ and $\ext^n_\mathcal{Y}$
simultaneously. We list some of these properties which will be used frequently in the sequel:
\begin{enumerate}
\item If $\frak B$ is admissible, then there is a natural isomorphism $\Ext^0_{\frak B} (-,-) \iso \hom_\Lambda (-,-)$,
and the natural comparison $\mathbb{Z}$-homomorphism $\eta^1_{M,N} : \Ext^1_{\frak B} (M,N) \arrow \Ext^1_{\Lambda } (M,N)$ is injective
for all $\Lambda$-modules $M$ and $N$ (see e.g.~\cite[Proposition 2.1]{hom.flat}). In particular, $\Ext^1_{\Lambda } (M,N)=o$ implies $\Ext^1_{\frak B} (M,N)=o$.

\item Let $n \geq 0$ be an integer. If $M$ and $N$ are $\Lambda$-modules, then
$\Ext^{\geqslant n+1}_{\frak B} (M,N)=o $
provided that $\xpd_\Lambda (M) \leq n$ or $\yid_\Lambda (N) \leq n$; this follows easily from the definition of $\Ext^n_{\frak B} (M,N)$.
Conversely, the equality $\Ext^{\geqslant n+1}_{\frak B} (M,-)=o$ implies $\xpd_\Lambda (M) \leq n$, and the equality
$\Ext^{\geqslant n+1}_{\frak B} (-,N)=o$ implies $\yid_\Lambda (N) \leq n$; see e.g.~\cite[Theorem~(3.9)]{holm.rel.Ext}.

\item If $\frak B$ is admissible, then $\Ext^n_{\frak B} (P,-)=\Ext^n_{\frak B} (-,I)=o$ for every projective $\Lambda$-module $P$ and injective
$\Lambda$-module $I$ because the functors $\hom_\Lambda (P,-)$ and $\hom_\Lambda (-,I)$ are exact.

\item Given a short $\frak B$-exact sequence
\[\xymatrix{ o \ar[r] & M' \ar[r] & M \ar[r] & M'' \ar[r] & o } \: , \]
of $\Lambda$-modules, there are long exact sequences
\[
\xymatrix{
o \ar[r] & \Ext^0_{\frak B} (M'' , - ) \ar[r] & \Ext^0_{\frak B} (M,-) \ar[r] & \Ext^0_{\frak B} (M',-) \ar `r[d]  `[l]  `[dlll] `[ll] [lld] \\
& \Ext^1_{\frak B} (M'',-) \ar[r] & \Ext^1_{\frak B} (M,-) \ar[r] & \Ext^1_{\frak B} (M',-) \ar[r] & \cdots
}\] 
and
\[
\xymatrix{
o \ar[r] & \Ext^0_{\frak B} (-,M') \ar[r] & \Ext^0_{\frak B} (-,M) \ar[r] & \Ext^0_{\frak B} (-,M'') \ar `r[d]  `[l]  `[dlll] `[ll] [lld] \\
& \Ext^1_{\frak B} (-,M') \ar[r] & \Ext^1_{\frak B} (-,M) \ar[r] & \Ext^1_{\frak B} (-,M'') \ar[r] & \cdots 
}\]
of $\mathbb{Z}$-modules; see e.g.~\cite[Theorems 8.2.3 and~8.2.5]{enochs.book}.
If, furthermore,  $\frak B$ is admissible, we may identify $\Ext^0_{\frak B} (-,-)$ with $\hom_\Lambda (-,-)$ in the above long exact sequences; cf. part~(i).

\item For every $\Lambda$-module $M$ and every integer $n \geq 0$, the functor $\Ext^n_{\frak B} (M,-)$ commutes with 
products and the functor $\Ext^n_{\frak B} (-,M)$ turns co-products into products. If, furthermore, $M$ has a degreewise finitely generated
left $\mathcal{X}$-resolution, then $\hom_\Lambda (M,-)$ and hence $\Ext^n_{\frak B} (M,-)$ commutes with direct sums.
\end{enumerate}
\eobs

It is worth recording the following straightforward consequence of~\ref{obs:BExt}-(iv) because of its frequent use later in the arguments.

\bdes{Dimension Shifting Lemma} \label{lem:Dim.shifting}
\em
Let $M$ be a $\Lambda$-module, $\frak B$ be an admissible $\hom$-balanced pair in $\Mod (\Lambda)$, and 
\[\xymatrix{ o \ar[r] & A \ar[r]^-{d_1} & B_1 \ar[r] &  \cdots \ar[r]^-{d_n} & B_n \ar[r] & C \ar[r] & o }\]
be a $\frak B$-exact sequence of $\Lambda$-modules. 
\begin{enumerate}
\item If $\Ext^{\geqslant 1}_{\frak B} (M,B_i)=o$ for all $1 \leq i \leq n$, then $\Ext^{n + i}_{\frak B} (M,A) \cong_\mathbb{Z} \Ext^i_{\frak B} (M,C)$
for all $i \geq 1$. 

\item If $\Ext^{\geqslant 1}_{\frak B} (B_i , M)=o$ for all $1 \leq i \leq n$, then $\Ext^{n + i}_{\frak B} (C,M) \cong_\mathbb{Z} \Ext^i_{\frak B} (A,M)$ for all $i \geq 1$.
\end{enumerate}
\edes

By part~(i) of~\ref{obs:BExt}, the functor $\Ext_\mathcal{B}^1$ associated with an admissible $\hom$-balanced pair $\mathfrak{B}$ is a subfunctors of $\Ext^1_\Lambda$. As it is explained below, such subfunctors of $\Ext^1_\Lambda$ correspond to the ``additive subfunctors of $\Ext^1_\Lambda$ with enough projectives and injectives'' in the sense of Auslander and Solberg~\cite[Section~1]{ASI}.

\bdes{Additive Subfunctors of $\mathbf{Ext}$ and Balanced Pairs} \label{F-exact}
Let $F$ be an additive subfunctor of the functor $\Ext^1_\Lambda (-,-)$, and for $\Lambda$-modules $A$ and $C$ identify $\Ext^1_\Lambda (C,A)$ with the abelian group of extensions of $A$ by $C$. Following~\cite[Section~1]{ASI}, a short exact sequence of $\Lambda$-modules of the form 
\begin{equation} \label{eq:F-exact}
\xymatrix{o \ar[r] & A \ar[r] & B \ar[r] & C \ar[r] & o}
\end{equation}
is said to be \emph{$F$-exact} if it belongs to $F (C,A)$. A $\Lambda$-module $P$ is called \emph{$F$-projective} if for any $F$-exact sequence~\eqref{eq:F-exact},
the induced sequence
\[\xymatrix{o \ar[r] & \hom_\Lambda (P,A) \ar[r] & \hom_\Lambda (P,B) \ar[r] & \hom_\Lambda (P,C) \ar[r] & o }\] 
is exact. The class of all $F$-projective $\Lambda$-modules is denoted by $\mathscr{P} (F)$. The notion of an \emph{$F$-injective module} is defined dually, and the class of all $F$-injective $\Lambda$-modules is denoted by $\mathscr{I} (F)$. Note that $\Proj_0 \subseteq \mathscr{P} (F)$ and $\mathcal{I}_0 \subseteq \mathscr{I} (F)$.

An additive subfunctor $F$ of $\Ext^1_\Lambda (-,-)$ is said to have \emph{enough projective} if for every $\Lambda$-module $M$ there exists an $F$-exact
sequence of the form $\xymatrixcolsep{0.5cm}\xymatrix{o \ar[r] & K \ar[r] & P \ar[r] & M \ar[r] & o}$ with $P \in \mathscr{P} (F)$. Dually, $F$ is said to have 
\emph{enough injectives} if for every $\Lambda$-module $M$ there exists an $F$-exact sequence of the
form $\xymatrixcolsep{0.5cm}\xymatrix{o \ar[r] & M \ar[r] & I \ar[r] & C \ar[r] & o}$ with $I \in \mathscr{I} (F)$.
It follows easily from the definition that if $F$ has enough projectives and enough injectives, then $\mathfrak{B} (F) : = \big ( \mathscr{P} (F) , \mathscr{I} (F) \big )$ is an admissible $\hom$-balanced pair, which is called the \emph{$\hom$-balanced pair associated with $F$}. It then follows from~\cite[Theorem~1.15]{ASI}
that the assignments $F \mapsto \mathfrak{B} (F)$ and $\mathfrak{B} \mapsto \Ext^1_{\frak B} (-,-)$ induce mutually inverse bijections between
``additive subfunctors $F$ of $\Ext^1_\Lambda (-,-)$ with enough projectives and injectives'' and ``admissible $\hom$-balanced pairs $\mathfrak{B}$ 
in $\Mod (\Lambda)$''.
\edes

In this paper we shall be particularly interested in the``Gorenstein $\hom$-balanced pair'' which is the Gorenstein analogue of the
standard $\hom$-balanced pair $(\mathcal{P}_0 , \mathcal{I}_0)$. We need to fix some notations before we formally introduce the Gorenstein $\hom$-balanced pair.
For unexplained notions from Gorenstein homological algebra we refer the reader to~\cite[Chapter 10]{enochs.book} and~\cite{chen.gorenstein}.

\bnotation \label{notation:GP} 
For an integer $n \geq 0$ we denote by $\GP_n$ the class of all $\Lambda$-modules of Gorenstein projective dimension at most $n$.
Likewise, we denote the class of  all $\Lambda$-modules of Gorenstein injective dimension at most $n$ by
$\GI_n$. The Gorenstein projective dimension of a $\Lambda$-module is denoted by $\Gpd_\Lambda (-)$.
\enotation

\bdes{Gorenstein $\mathbf{Hom}$-balanced Pair}
\label{df:vGorestein} 
\label{df:Gbalanced} 
\label{Gext}
A ring $\Lambda$ is called \emph{virtually Gorenstein} if $\mathfrak{G}:=(\GP_0 , \GI_0)$ is a (necessarily admissible) $\hom$-balanced pair.
Such rings were first introduced and studied in the context of representation theory of artin algebras by Beligiannis and Reiten~\cite{BR,beligiannis.CM}.
If $\Lambda$ is a virtually Gorenstein ring, the $\hom$-balanced pair $\mathfrak{G}$ is referred to as the \emph{Gorenstein $\hom$-balanced pair} 
or the \emph{$G$-balanced pair} for short, and  $\mathfrak{G}$-exact sequences in the sense of~\ref{df:B-exactness} are simply called \emph{$G$-exact}. 
Note that $G$-exact complexes are necessarily exact. Following~\cite{enochs.Gext} and~\cite{holm.Gext}, the relative $\Ext$-functors associated with $\frak G$ will be denoted by $\Gext^n_\Lambda (-,-)$ for all $n \in \mathbb{Z}$.
\edes

\bdes{Ubiquity of Virtually Gorenstein Rings} \label{CM.finite}
A classical example of a virtually Gorenstein ring is an Iwanaga-Gorenstein ring  (see e.g.~\cite[Theorem 2.1]{enochs.Gext} and~\cite[Theorem 12.1.4]{enochs.book}),
and this indeed justifies the nomenclature ``virtually Gorenstein'' to some extent. virtually Gorenstein artin algebras are abundant. Indeed, any artin algebra of finite representation type is virtually Gorenstein (see~\cite[Example 8.4.]{beligiannis.CM}), and by~\cite[Theorem~8.11 and Theorem~8.12]{beligiannis.CM}
the property of ``being virtually Gorenstein'' for artin algebras is preserved under derived equivalences and stable equivalences of Morita type.
See~\cite{belig.krause.thick} for an example of an artin algebra which is \emph{not} virtually Gorenstein.

Virtually Gorenstein algebras of ``finite CM-type'' are of particular importance in this work. Recall that an artin algebra is said to be of \emph{finite Cohen-Macaulay type}, or of \emph{finite CM-type} for short, if the class $\Gp_0$ is of finite representation type, i.e. there are only finitely many indecomposable finitely generated Gorenstein projective modules up to isomorphism.
\edes 

We conclude this section with the following lemma about $G$-exactness of direct limits of short $G$-exact sequences.

\blem \label{lem:dirlim.of.Gexact}
Let $\Lambda$ be a virtually Gorenstein ring.
\begin{enumerate}
\item The direct sum of any family of short $G$-exact sequences of $\Lambda$-modules is $G$-exact.

\item If every Gorenstein projective $\Lambda$-module decomposes as a direct sum of finitely generated Gorenstein projective modules (see~Remark~\ref{rem:GP.over.CM}), then the direct limit of any directed system of short $G$-exact sequences is $G$-exact.
\end{enumerate}
\elem

\bproof
Part~(i): Let $\big \{ \xymatrix{ \delta_i : = o \ar[r] & A_i \ar[r]^-{f_i} & B_i \ar[r]^-{g_i} & C_i \ar[r] & o } \big \}_{i \in I}$ be a family of 
short $G$-exact sequences of $\Lambda$-modules and denote by $\coprod_{i \in I} \delta_i$ the direct sum of this family. Note that 
$\coprod_{i \in I} \delta_i$ is exact. To check that it is also $G$-exact observe that for every Gorenstein injective $\Lambda$-module $I$ 
we have the canonical isomorphism
\[ \hom_\Lambda \Big ( \coprod_{i \in I} \delta_i , I \Big ) \cong \prod_{i \in I} \hom_\Lambda ( \delta_i , I ) \]
of $\mathbb{Z}$-complexes. Since each $\delta_i$ is $G$-exact, each complex $\hom_\Lambda ( \delta_i , I )$ is exact and hence the 
product $\prod_{i \in I} \hom_\Lambda  ( \delta_i , I  ) $ is exact. Consequently $\coprod_{i \in I} \delta_i$ is $G$-exact.

Part~(ii): Let $\big \{ \xymatrix{ \delta_i : = o \ar[r] & A_i \ar[r]^-{f_i} & B_i \ar[r]^-{g_i} & C_i \ar[r] & o } \big \}_{i \in I}$ be a direct system of
short $G$-exact sequences of $\Lambda$-modules and let $\delta := \dirlim_{i \in I} \delta_i$. It is clear that $\delta$ is exact because direct limit is 
an exact functor on the category of modules. Now let $G$ be a Gorenstein projective $\Lambda$-module. By the hypothesis,
$G = \coprod_{j \in J} G_j$ where each $G_j$ is finitely generated Gorenstein projective for all $j \in J$. Thus, $\hom_\Lambda (G,-) \cong \prod_{j \in J} \hom_\Lambda (G_j , - )$, and since each $G_i$ is finitely generated, by well-known canonical isomorphisms in module theory, we have the following isomorphisms of complexes of abelian groups:

\begin{align*}
\hom_\Lambda (G, \delta) & \cong \prod_{j \in J} \hom_\Lambda (G_j , \delta ) \\
& \cong \prod_{j \in J} \dirlim_{i \in I} \hom_\Lambda (G_j , \delta_i) \: .
\end{align*}
For every $i \in I$ and $j \in J$, the complex $\hom_\Lambda (G_j , \delta_i)$ is exact because $\delta_i$ is $G$-exact by the hypothesis. Therefore, for every $j \in J$ the sequence $\dirlim_{i \in I} \hom_\Lambda (G_j , \delta_i)$ is exact
being the direct limit of exact sequences of abelian groups. Hence, $\prod_{j \in J} \dirlim_{i \in I} \hom_\Lambda (G_j , \delta_i)$ is exact being a product of exact sequences of abelian groups. It then follows that $\hom_\Lambda (G, \delta)$ is exact and this complete the proof of $G$-exactness of $\delta$.
\eproof

\brem \label{rem:GP.over.CM}
By~\cite[Theorem~4.10]{beligiannis.CMtype}, every Gorenstein projective module over a virtually Gorenstein CM-finite artin algebra decomposes into a direct sum of finitely generated Gorenstein projective modules. Hence the assumption of part~(ii) of Lemma~\ref{lem:dirlim.of.Gexact} is satisfied over CM-finite virtually Gorenstein artin algebras.
\erem

\section{Relative Orthogonal Classes and Eklof Lemma}
\label{sec:orth}

In this section we discuss relative $\Ext$-orthogonal classes, which are classes of modules defined as roots of relative Ext-functors associated with
$\hom$-balanced pairs. The key result of this section, namely Lemma~\ref{Eklof.lemma.Rel}, is the relative version of Eklof Lemma (cf.~\cite[Lemma 6.2]{GT})
which guarantees closure of relative left $\Ext$-orthogonal classes under certain filtrations; cf.~\ref{df:rel.filtration}.

\bconv
Throughout this section, $\frak B := (\mathcal{X} , \mathcal{Y})$ denotes an \emph{admissible $\hom$-balanced pair} in the sense defined in~\ref{df:admissible.Bpair}. 
\econv

\bdes{Relative $\mathbf{Ext}$-orthogonal Classes} \label{note:*perp}
For a class $\mathcal{C}$ of $\Lambda$-modules, let
\[ \mathcal{C}^{\bperp} : = \left \{ M \in \Mod (\Lambda) : \text{$\Ext^1_{\frak B} (C , M) = o$ for all $C \in \mathcal{C}$} \right \}  \]
and
\[ {}^{\bperp} \mathcal{C} : = \left \{ M \in \Mod (\Lambda) : \text{$\Ext^1_{\frak B} (M,C) = o$ for all $C \in \mathcal{C}$} \right \} \: . \]
These are called the \emph{right} and the \emph{left} \emph{$\Ext_{\mathfrak{B}}$-orthogonal classes} of $\mathcal{C}$ respectively.
In the special case where $\mathcal{C}$ consists  of only one element $T$, we let
$T^\bperp := \{ T \}^\bperp$ and ${}^\bperp T :={}^\bperp \{ T \}$ for convenience. 
\edes

\bdes{Notational Remark}
Two special cases of $\mathfrak{B}$ are of particular interest with regard to relative orthogonal classes:
When $\mathfrak{B}$ is the standard $\hom$-balanced pair~\ref{standard.balanced.pair}, the left/right $\Ext_{\mathfrak{B}}$-orthogonal classes become the standard left/right $\Ext$-orthogonal classes defined in terms of vanishing of $\Ext^1$. In this case, we simply replace $\perp\hspace{-0.2cm}\mathfrak{B}$ with $\perp$, and for example we write $\mathcal{C}^{\perp}$ instead of $\mathcal{C}^{\bperp}$. Likewise, when $\mathfrak{B}$ is the Gorenstein $\hom$-balanced pair from~\ref{df:Gbalanced}, 
we replace $\perp\hspace{-0.2cm}\mathfrak{B}$ with $\mathrm{G}\hspace{-0.2cm}\perp$ in the notation and write for example $\mathcal{C}^{\Gperp}$ instead of $\mathcal{C}^{\bperp}$.
\edes

\bobs \label{obs:*perp}
For a class $\mathcal{C}$ of $\Lambda$-modules it follows readily from basic properties of $\Ext_{\frak B}$-functor
mentioned in~\ref{obs:BExt} that:
\begin{enumerate}
\item $\mathcal{C}^\bperp$ is closed under $\frak B$-extensions, direct summands, and products. Furthermore, $\mathcal{C}^{\bperp}$ contains $\mathcal{I}_0$
and $\mathcal{Y}$.

\item ${}^\bperp \mathcal{C}$ is closed under  $\frak B$-extensions, direct summands, and coproducts. Furthermore, ${}^{\bperp} \mathcal{C}$ 
contains $\mathcal{P}_0$ and $\mathcal{X}$.

\item If $\mathcal{C}$ is a \emph{set} (as opposed to a class), then we can set $\mathcal{C}^\bperp = M^\bperp$ and
${}^\bperp \mathcal{C}={}^\bperp N$  where $M=\coprod_{C \in \mathcal{C}} C$ and $N = \prod_{C \in \mathcal{C}} C$; cf. part~(v) of~\ref{obs:BExt}.

\item $\mathcal{C}^{\perp} \subseteq \mathcal{C}^\bperp$ and $ {}^\perp \mathcal{C} \subseteq {}^\bperp \mathcal{C}$; cf. part~(i) of~\ref{obs:BExt}.

\item The inclusions $\mathcal{C} \subseteq {}^\bperp (\mathcal{C}^\bperp)$ and $\mathcal{C} \subseteq ({}^\bperp \mathcal{C})^\bperp$ always hold.
\end{enumerate}
Furthermore, if $\mathcal{D}$ is another class of $\Lambda$-modules with $\mathcal{C} \subseteq \mathcal{D}$, 
then $\mathcal{D}^\bperp \subseteq \mathcal{C}^\bperp$ and ${}^\bperp \mathcal{D} \subseteq {}^\bperp \mathcal{C}$. 
From this and part~(v), one can easily deduce the equalities
\begin{equation} \label{eq:*perp.equality}
{}^\bperp \big ( ({}^\bperp \mathcal{C})^\bperp \big ) = {}^\bperp \mathcal{C} \quad , \quad 
\big ( {}^\bperp (\mathcal{C}^\bperp) \big )^\bperp = \mathcal{C}^\bperp
\end{equation}
for every class $\mathcal{C}$ of $\Lambda$-modules.
\eobs

A well-known lemma due to Eklof~\cite{eklof}, states that left $\Ext$-orthogonal classes 
of modules  are closed under filtrations (see~\cite[Lemma 6.2]{GT} for a precise statement).
We are next going to prove a relative analogue of this result---see~\ref{Eklof.lemma.Rel} below---which states that left $\Ext_\mathfrak{B}$-orthogonal classes are closed under ``$\frak B$-proper filtrations''  in the sense defined below.

\bdes{$\mathfrak{B}$-proper Filtrations} \label{df:rel.filtration}
\label{G.proper.filt}
Recall that a \emph{continuous chain} of $\Lambda$-modules of length $\kappa$ ($\kappa$ is an ordinal)
is a family $\{ M_\alpha \}_{\alpha \leq \kappa}$ of $\Lambda$-modules such that $M_0 = o$, $M_{\alpha} \subseteq M_{\alpha +1}$ for any 
ordinal $\alpha < \kappa$, and $M_\alpha = \bigcup_{\beta < \alpha} M_\alpha$ for any limit ordinal $\alpha \leq \kappa$. Such a continuous chain is 
said to be \emph{$\mathfrak{B}$-proper} if for every ordinal $\alpha < \kappa$, the short exact sequence
\begin{equation} \label{eq:continuous.chain}
\xymatrix{ o \ar[r] & M_\alpha \ar@{^{(}->}[r] & M_{\alpha+1} \ar[r] & M_{\alpha+1} / M_\alpha \ar[r] & o }
\end{equation}
is $\frak B$-exact. Given a class $\mathcal{C}$ of $\Lambda$-modules, a $\Lambda$-module $M$ is said to be \emph{$\frak B$-properly $\mathcal{C}$-filtered},
or to be a \emph{transfinite $\frak B$-extension of modules in $\mathcal{C}$} if there is a $\frak B$-proper continuous chain $\{ M_{\alpha} \}_{\alpha \leq \kappa}$ 
of $\Lambda$-submodules of $M$ such that $M_{\kappa} = M$ and for every ordinal $\alpha < \kappa$, 
the $\Lambda$-module $M_{\alpha+1}/M_\alpha$ is isomorphic to some member of $\mathcal{C}$. The family $\{ M_{\alpha} \}_{\alpha \leq \kappa}$ is then 
called a \emph{$\frak B$-proper $\mathcal{C}$-filtration} of $M$. In the special case where $\kappa$ is a finite ordinal (i.e. a natural number),
$M$ is said to be \emph{finitely $\frak B$-properly $\mathcal{C}$-filtered}.
\edes

\brem
When $\mathcal{B}$ is the $G$-balanced pair over a virtually Gorenstein ring, we replace the qualifier
``$\mathfrak{B}$-properly'' with ``$G$-properly''. Thus, for example, we speak of ``$G$-proper $\mathcal{C}$-filtrations'' instead of
``$\mathfrak{B}$-proper $\mathcal{C}$-filtrations" when $\mathfrak{B}=(\GP_0,\GI_0)$.
\erem

\bdes{Relative Eklof Lemma} \label{Eklof.lemma.Rel}
\em
Let $\mathcal{C}$ be a class of $\Lambda$-modules. If a $\Lambda$-module $M$ is $\mathfrak{B}$-properly ${}^\bperp \mathcal{C}$-filtered,
then $M \in {}^{\bperp}\mathcal{C}$.
\edes

The proof of the lemma is essentially along the lines of the proof of the ``absolute version'' of Eklof Lemma (see e.g.~\cite[Lemma 6.2]{GT}). 
However, we provide a detailed proof for reader's convenience.

\bproof
For any $\Lambda$-module $M$, $M \in {}^\bperp \mathcal{C}$ if and only if $M \in {}^\bperp C$ for all $C \in \mathcal{C}$.
Furthermore, since ${}^\bperp \mathcal{C} \subseteq {}^\bperp C$ for all $C \in \mathcal{C}$, every $\frak B$-properly ${}^\bperp \mathcal{C}$-filtered 
module is also $\frak B$-properly ${}^\bperp C$-filtered for all $C \in \mathcal{C}$. Consequently, 
in order to prove the assertion it suffices to show that ${}^\bperp C$ is closed under transfinite $\frak B$-extensions for any $\Lambda$-module $C$.

Let $\left \{ M_\alpha \right \}_{\alpha \leq \kappa}$ be a $\frak B$-proper ${}^\bperp C$-filtration of $M$. We use transfinite induction to show that for 
every ordinal $\alpha \leq \kappa$, $M_\alpha \in {}^\bperp C$, which in particular implies $M=M_\kappa \in {}^\bperp C$. 

Let $\alpha$ be an ordinal such that $\alpha < \kappa$. If $\alpha=0$, then $M_\alpha = o$ is certainly in ${}^\bperp C$. 
Furthermore, if $M_\alpha$ belongs to ${}^\bperp C$, then $M_{\alpha +1}$ also belongs to ${}^\bperp C$ because ${}^\bperp C$ is closed 
under $\frak B$-extensions and the short exact sequence
\[\xymatrix{ o \ar[r] & M_\alpha \ar@{^{(}->}[r] & M_{\alpha+1} \ar[r] & M_{\alpha+1} / M_{\alpha} \ar[r] & o }\]
is a $\frak B$-extension, whose end terms belong to ${}^\bperp C$. Thus, to finish the proof we need to show that if $\alpha \leq \kappa$ is a limit ordinal 
and $M_\beta \in {}^\bperp C$ for all ordinals $\beta < \alpha$, then $M_\alpha \in {}^\bperp C$.
Let $j : C \arrow Y$ be a $\mathcal{Y}$-preenvelope of $C$ and note that the associated $\frak B$-exact sequence
$\xymatrix{o \ar[r] & C \ar[r]^-j & Y \ar[r]^-p & W \ar[r] & o}$ induces the exact sequence
\[\xymatrix{ \hom_\Lambda (-,Y) \ar[r]^-{p_\ast} & \hom_\Lambda (-,W) \ar[r] & \Ext^1_{\frak B} (-,C) \ar[r] & o } \]
by part~(iv) of Observation~\ref{obs:BExt}. Therefore, $M_\alpha$ belongs to ${}^\bperp C$ if and only if
\[\xymatrix{ \hom_\Lambda (M_\alpha,Y) \ar[r]^-{p_\ast} & \hom_\Lambda (M_\alpha,W) }\]
is surjective. Thus, the proof of $M_\alpha \in {}^\bperp C$ reduces to the following ``lifting problem'':
\emph{Given a $\Lambda$-homomorphism $f : M_\alpha \arrow W$, there 
exists a $\Lambda$-homomorphism $g : M_\alpha \arrow Y$ making the diagram}
\[\xymatrix{
& M_\alpha \ar[d]^-f \ar@{-->}[dl]_-g \\
Y \ar@{->>}[r]_-p & W
}\]
\emph{commute}. For every ordinal $\beta \leq \alpha$ let $f_\beta :=f |_{M_\beta}$. Since $M_\alpha = \bigcup_{\beta < \alpha} M_\beta$,
in order to construct $g$ it suffices to construct a family $\left \{ g_\beta : M_\beta \arrow Y \right \}_{\beta < \alpha}$ of $\Lambda$-homomorphisms 
with the following compatibility properties:
\begin{enumerate}
\item[(1)] $p \circ g_\beta = f_\beta$ and

\item[(2)] for all ordinals $\gamma < \beta < \alpha$, ${g_\beta}|_{M_\gamma}=g_\gamma$.
\end{enumerate}
The family  $\left \{ g_\beta : M_\beta \arrow Y \right \}_{\beta < \alpha}$ can be constructed by transfinite induction on $\beta < \alpha$ as follows:
For $\beta = 0$ simply let $g_0 = 0$. Assume now that $g_{\beta} : M_\beta \arrow Y$ is already constructed for some ordinal $\beta < \alpha$ and
note that the short $\mathfrak{B}$-exact sequences
\[\xymatrix{ o \ar[r] & C \ar[r]^-j & Y \ar[r]^-p  & W \ar[r]& o }\]
and
\[\xymatrix{ o \ar[r] & M_\beta \ar@{^{(}->}[r] & M_{\beta+1} \ar[r] & \frac{M_{\beta+1}}{M_\beta} \ar[r] & o }\]
induce the commutative diagram 
\[
\xymatrixcolsep{0.65cm}
\xymatrixrowsep{0.65cm}
\xymatrix{
& o \ar[d] & o \ar[d] & o \ar[d] & \\ 
o \ar[r] & \hom_\Lambda \big (\frac{M_{\beta+1}}{M_\beta} , C \big ) \ar[r]^-{j_1} \ar[d] & \hom_\Lambda \big ( \frac{M_{\beta+1}}{M_\beta} , Y \big ) \ar[r]^-{p_1} \ar[d] & 
\hom_\Lambda \big ( \frac{M_{\beta+1}}{M_\beta} , W \big ) \ar[r] \ar[d] & o \\
o \ar[r] & \hom_\Lambda (M_{\beta+1} , C) \ar[r]^-{j_2} \ar[d]^-{i_1} & \hom_\Lambda (M_{\beta+1} , Y) \ar[r]^-{p_2} \ar[d]^-{i_2} & 
\hom_\Lambda (M_{\beta+1} , W)  \ar[d]^-{i_3} &  \\
o \ar[r] & \hom_\Lambda (M_\beta , C) \ar[r]^-{j_3} \ar[d] & \hom_\Lambda (M_\beta , Y) \ar[r]^-{p_3} \ar[d] & \hom_\Lambda (M_\beta , W)  &  \\
          & o                                                & o                                               &                                                 &          
}\]
with exact rows and columns in view of $\Ext^1_\mathfrak{B} (-,Y)=o$ and $M_{\beta+1}/M_\beta \in {}^{\bperp} C$.
Then the equalities ${f_{\beta+1}}|_{M_\beta} = p \circ g_\beta =f_\beta$ translate into $i_3 (f_{\beta+1})=p_3 (g_\beta)=f_\beta$, and 
an easy diagram chasing shows that there exists $g_{\beta+1} \in \hom_\Lambda (M_{\beta+1} , Y)$ such that $p_2 (g_{\beta+1})=f_{\beta+1}$ 
and $i_2 (g_{\beta+1})=g_{\beta}$. These identities can be rephrased as $p \circ g_{\beta+1}=f_{\beta+1}$ and $g_{\beta+1} |_{M_\beta}=g_{\beta}$, which are the desired conditions (1) and (2) mentioned above. Finally, if $\beta < \alpha$ is a limit ordinal and $g_\gamma : M_\gamma \arrow I$ is constructed as required for all $\gamma < \beta$, then we may let $g_\beta = \bigcup_{\gamma < \beta} g_\gamma$. The family $\{ g_\beta : M_\beta \arrow I \}_{\beta < \alpha}$ now has the desired 
properties and this finishes the proof.
\eproof
 
\section{Relative Cotorsion Pairs}
\label{sec:Rel.cotor}

In this section we consider cotorsion pairs relative to $\hom$-balanced pairs. Cotorsion pairs---first introduced by Salce~\cite{salce} as an analogue of (non-hereditary) torsion pairs---play an important role in approximation theory of modules due to their close connection with special approximations.
It was proved by Salce~\cite{salce} that the components of a cotorsion pair are ``dual'' in the sense that the left-hand side component is special precovering if and only if the right-hand side component is special preenveloping; cf.~\ref{salce.lemma}. Later on, it was proved by Eklof and Trlifaj~\cite{trlifaj.Extvanish} that every cotorsion pair cogenerated by a set of modules is complete; cf.~\ref{thm:eklof.trlifaj*}. These two facts, among other things, indicate that cotorsion pairs are an invaluable machinery for detecting or producing classes of modules which provide special approximations. The reader is referred to~\cite{GT} for a detailed discussion on cotorsion pairs.

Given a $\hom$-balanced pair $\mathfrak{B}$ in $\Mod (\Lambda)$, one can use $\Ext^1_{\mathfrak{B}}$, instead of $\Ext^1_\Lambda$ in the definition of a cotorsion pair, to define ``relative cotorsion pairs'' with respect to $\frak B$; see Definition~\ref{df:cotorsion.pair}. 
Cotorsion pairs relative to a $\hom$-balanced pair have already been considered in~\cite{rel.cotorsion} and it was shown therein that
some basic properties of cotorsion pairs directly carry over into the relative setting; see e.g.~\cite[Section~3]{rel.cotorsion}. 
The aim of this section is to promote the analogy between absolute and relative cotorsion pairs by establishing the relative counterpart of 
Eklof-Trlifaj Completeness Theorem~\cite[Theorem 2]{trlifaj.Extvanish}. The results of this section will then be used in the subsequent sections to study
``Gorenstein tilting approximations''. We begin with the definition of relative cotorsion pairs and recalling some of their basic properties from~\cite{rel.cotorsion}.

\bdes{Relative Cotorsion Pairs} 
\label{df:cotorsion.pair} 
\label{df:rel.cotorsionpair.gen}
Let $\frak B : = (\mathcal{X} , \mathcal{Y})$ be a $\hom$-balanced pair in $\Mod (\Lambda)$.
A pair $\mathfrak{C} : = ( \mathcal{L} , \mathcal{R} )$ of full subcategories of $\Mod (\Lambda)$ is called a \emph{cotorsion pair relative to $\frak B$} or
simply a \emph{$\mathfrak{B}$-cotorsion pair} if $\mathcal{L}^\bperp =\mathcal{R}$ and $\mathcal{L}={}^\bperp \mathcal{R}$. In this case, the class $\mathcal{K}_{\mathfrak{C}}:=\mathcal{L} \cap \mathcal{R}$ is called the \emph{kernel} of the relative cotorsion pair $\mathfrak{C}$.
The class $\mathcal{L}$ is referred to as the \emph{left component} of $\mathfrak{C}$ and likewise the class $\mathcal{R}$ is referred to as the \emph{right component} of $\mathfrak{C}$. It is clear that $\big (\mathcal{X} , \Mod (\Lambda) \big )$ and $\big ( \Mod (\Lambda) , \mathcal{Y} \big )$ are always $\frak B$-cotorsion pairs. These are called  \emph{trivial $\frak B$-cotorsion pairs}, and they are indeed the relative counterparts of the trivial cotorsion
pairs $\big ( \mathcal{P}_0  , \Mod (\Lambda) \big )$ and $\big ( \Mod (\Lambda) , \mathcal{I}_0  \big )$.

Cotorsion pairs relative to $(\mathcal{P}_0 , \mathcal{I}_0)$ are simply the  cotorsion pairs in the usual sense. 
Cotorsion pairs relative to the $G$-balanced pair $(\GP_0 , \GI_0)$ are referred to as \emph{$G$-cotorsion pairs}.

Given a class $\mathcal{C}$ of $\Lambda$-modules, we may form by~\eqref{eq:*perp.equality}
the $\mathfrak{B}$-cotorsion pair $\big ( {}^\bperp (\mathcal{C}^\bperp) , \mathcal{C}^\bperp \big )$
\emph{cogenerated by $\mathcal{C}$} and the $\mathfrak{B}$-cotorsion pair $\big ( {}^\bperp \mathcal{C} , ({}^\bperp \mathcal{C})^\bperp \big )$  \emph{genetated by $\mathcal{C}$}. A $\mathfrak{B}$-cotorsion pair $(\mathcal{L} , \mathcal{R})$ is said to be \emph{cogenerated by $\mathcal{C}$} if $\mathcal{R}=\mathcal{C}^\bperp$ and it is said to be \emph{generated by $\mathcal{C}$} if ${}^\bperp \mathcal{C} = \mathcal{L}$. 
\edes

The following proposition about  basic closure properties of relative cotorsion pairs follows easily from the Definition~\ref{df:cotorsion.pair} and Observarion~\ref{obs:*perp}.

\bprop \label{prop:cotorsion.closure}
Let $\frak B  = (\mathcal{X} , \mathcal{Y})$ be a $\hom$-balanced pair in $\Mod (\Lambda)$.
If $\big ( \mathcal{L} , \mathcal{R} \big )$ is a $\mathfrak{B}$-cotorsion pair in $\Mod (\Lambda)$, then:
\begin{enumerate}
\item $\mathcal{L}$ contains $\mathcal{X}$ and it is closed under $\frak B$-extensions, direct sums and direct summands.

\item $\mathcal{R}$ contains $\mathcal{Y}$ and it is closed under $\frak B$-extensions, products, and direct summands.
\end{enumerate}
\eprop

\bconv
From now on, until the end of this section, we assume that $\frak B := (\mathcal{X} , \mathcal{Y})$ is an \emph{admissible $\hom$-balanced pair}
in the sense of~\ref{df:admissible.Bpair}. 
\econv

\bdes{$\frak B$-(co)resolving Subcategories} \label{df:*resolving}
A full subcategory $\mathcal{C}$ of $\Mod (\Lambda)$ is said to be
\emph{$\frak B$-resolving} if it contains $\mathcal{X}$ and given a short $\frak B$-exact sequence
\begin{equation} \label{eq:*resolving}
\xymatrix{ o \ar[r] & M' \ar[r] & M \ar[r] & M'' \ar[r] & o }
\end{equation}
in $\Mod (\Lambda)$ with $M'' \in \mathcal{C}$, $M' \in \mathcal{C}$ if and only if $M \in \mathcal{C}$. Dually, $\mathcal{C}$ is said to be
\emph{$\frak B$-coresolving} if it contains $\mathcal{Y}$ and given a short $\frak B$-exact sequence in $\Mod (\Lambda)$ as~\eqref{eq:*resolving}
with $M' \in \mathcal{C}$, $M \in \mathcal{C}$ if and only if $M'' \in \mathcal{C}$.
Clearly for $\mathfrak{B}=(\mathcal{P}_0 , \mathcal{I}_0)$ a $\mathfrak{B}$-(co)resolving class is a ``(co)resolving class'' in the usual sense.
\edes

\brem \label{rem:*syzygy.closed}
It is easy to see that $\frak B$-resolving subcategories are \emph{$\mathcal{X}$-syzygy closed} in the sense that they contain all the relative syzygy modules of 
left $\mathcal{X}$-resolutions of their members.
Dually, $\frak B$-coresolving subcategories are \emph{$\mathcal{Y}$-cosyzygy closed} in the sense  that they contain all the relative cosyzygy modules of right $\mathcal{Y}$-resolutions of their members.
\erem

\bobs \label{obs:*inf.classes}
For a class $\mathcal{C}$ of $\Lambda$-modules let
\begin{eqnarray*}
\mathcal{C}^{\bperp_\infty} & : = & \left \{ M \in \Mod (\Lambda) \mid \hbox{$\Ext^{\geqslant 1}_\frak B (C , M) = o$ for all $C \in \mathcal{C}$}   \right \} \: , \\
{}^{\bperp_\infty} \mathcal{C} & : = & \left \{ M \in \Mod (\Lambda) \mid  \hbox{$\Ext^{\geqslant 1}_\frak B (M,C) = o$ for all $C \in \mathcal{C}$}   \right \} \: .
\end{eqnarray*}
It is easy to see that:
\begin{enumerate}
\item  the class $\mathcal{C}^{\bperp_\infty}$ is $\frak B$-coresolving and hence $\mathcal{Y}$-cosyzygy closed;

\item the class ${}^{\bperp_\infty} \mathcal{C}$ is $\frak B$-resolving and hence $\mathcal{X}$-syzygy closed;

\item $\mathcal{C}^{\bperp_\infty} \subseteq \mathcal{C}^\bperp$ and the equality holds provided that $\mathcal{C}$ is $\mathcal{X}$-syzygy closed;

\item ${}^{\bperp_\infty} \mathcal{C} \subseteq {}^\bperp \mathcal{C}$ and the equality holds provided that $\mathcal{C}$ is $\mathcal{Y}$-cosyzygy closed.
\end{enumerate}
\eobs

The following proposition states that there is a duality between $\frak B$-resolving and $\frak B$-coresolving properties of the 
components in a $\mathfrak{B}$-cotorsion pair.

\bprop[{\cite[Theorem~3.8]{rel.cotorsion}}] \label{prop:*rogas}
The following statements about a $\mathfrak{B}$-cotorsion pair $(\mathcal{L} , \mathcal{R})$ are equivalent :
\begin{enumerate}
\item $\mathcal{L}$ is $\frak B$-resolving;

\item $\mathcal{R}$ is $\frak B$-coresolving;

\item $\Ext^{\geqslant 1}_\frak B (A,B) = o$ for all $A \in \mathcal{L}$ and $B \in \mathcal{R}$.
\end{enumerate}
\eprop

\bdes{Hereditary $\frak B$-cotorsion Pairs} \label{df:hereditary.cotorsion}
A $\frak B$-cotorsion pair $(\mathcal{L} , \mathcal{R})$ in $\Mod (\Lambda)$ is said to be \emph{hereditary}
if it satisfies one of the equivalent conditions (i)--(iii) in~\ref{prop:*rogas}.
\edes

\bdes{$\frak B$-special Approximations} \label{special.approximation}
\label{df:G.special.approx}
Let $\mathcal{C}$ be a class of $\Lambda$-modules. A $\Lambda$-epimorphism $f : C \arrow M$ with $C \in \mathcal{C}$ is said to be a 
\emph{$\frak B$-special $\mathcal{C}$-precover} of $M$ if the sequence
\[\xymatrix{ o \ar[r] & \ker (f) \ar[r] & C \ar[r]^-f & M \ar[r] & o }\]
is $\frak B$-exact and $\ker (f) \in \mathcal{C}^{\bperp}$. Note that in this case $f : C \arrow M$ is a $\mathcal{C}$-precover. 
Dually, a $\Lambda$-monomorphism $f : M \arrow C$ with $C \in \mathcal{C}$ is said to be a
\emph{$\frak B$-special $\mathcal{C}$-preenvelop} of $M$ if the sequence
\[\xymatrix{ o \ar[r] & M \ar[r]^-{f} & C \ar[r] & \coker (f) \ar[r] & o }\]
is $\frak B$-exact and $\coker (f) \in {}^{\bperp} \mathcal{C}$. If every $\Lambda$-module has a $\frak B$-special $\mathcal{C}$-precover (respectively, 
a $\frak B$-special $\mathcal{C}$-preenvelope), then $\mathcal{C}$ is called a \emph{$\frak B$-special precovering class} (respectively, 
a \emph{$\frak B$-special preenveloping class}). 

We use the term ``$\frak B$-special approximation'' to refer to a $\frak B$-special precover or
a $\frak B$-special preenvelop. When $\mathfrak{B}$ is the standard $\hom$-balanced pair~\ref{standard.balanced.pair}, the $\frak B$-special approximations are simply called \emph{special approximations}. Furthermore, if $\mathfrak{B}$ is the $G$-balanced pair from~\ref{df:Gbalanced}, 
$\mathfrak{B}$-special approximations are called \emph{$G$-special approximations}.
\edes

As it was indicated at the beginning of this section, the importance of cotorsion pairs stems in large part from the
Salce Lemma (see~\cite{salce} and~\cite[Lemma~5.20]{GT}) and Eklof-Trlifaj Completeness Theorem (see~\cite{trlifaj.Extvanish} and~\cite[Theorem 6.11]{GT}). We are now going to discuss the relative counterparts of these results for $\frak B$-cotorsion pairs. We begin by recalling the relative version of Salce Lemma which has already appeared in~\cite[Theorem~3.11]{rel.cotorsion}.

\bdes{Relative Salce Lemma} \label{salce.lemma}
\em 
The following statements are equivalent for a $\frak B$-cotorsion pair $(\mathcal{L} , \mathcal{R})$:
\begin{enumerate}
\item $\mathcal{L}$ is $\frak B$-special precovering;

\item $\mathcal{R}$ is $\frak B$-special enveloping.
\end{enumerate}
\edes

\bdf \label{df:complete.cotorsion*}
A $\frak B$-cotorsion pair in $\Mod (\Lambda)$ satisfying one of the equivalent  conditions (i)  or (ii) in~\ref{salce.lemma} is 
called a \emph{complete $\frak B$-cotorsion pair}. 
\edf

We are now in a position to state and prove the relative version of Eklof-Trlifaj Completeness Theorem which guarantees
that $\frak B$-cotorsion pairs cogenerated by a \emph{set} of modules are complete.

\bdes{Relative Eklof-Trlifaj Completeness Theorem} \label{thm:eklof.trlifaj*}
\em 
Let $\mathcal{S}$ be a \emph{set} of $\Lambda$-modules and $\mathcal{S}^+:=\mathcal{S} \cup \mathcal{X}$.
\begin{enumerate}
\item For every $\Lambda$-module $M$ there is an exact short $\mathfrak{B}$-exact sequence
\[\xymatrix{o \ar[r] & M \ar[r]^-g & P  \ar[r] & N \ar[r] & o }\]
of $\Lambda$-homomorphisms wherein $P \in \mathcal{S}^\bperp$ and $N$ is $\frak B$-properly $\mathcal{S}$-filtered. 
In particular, $g : M \arrow P$ is a $\frak B$-special $\mathcal{S}^\bperp$-preenvelope of $M$.

\item The $\frak B$-cotorsion pair $\big ( {}^\bperp (\mathcal{S}^\bperp) , \mathcal{S}^\bperp \big )$ is complete
and the class  ${}^\bperp (\mathcal{S}^\bperp)$ consists precisely of all $\Lambda$-modules which are direct summands of 
$\mathfrak{B}$-properly $\mathcal{S}^+$-filtered modules.
\end{enumerate}
\edes

\bproof
Part~(i):~Let $S : = \coprod_{X \in \mathcal{S}} X$ and note that $\mathcal{S}^\bperp  = S^\bperp$ by part~(v) of Observation~\ref{obs:BExt}. 
Therefore, we may continue the argument assuming without loss of generality that $\mathcal{S}=\left \{ S \right \}$. 
We proceed in several steps to prove the assertion in part~(i) of the theorem:\\

\noindent \textsc{Step~1.} Since $\frak B$ is admissibe, there is an exact short $\frak B$-exact sequence
\begin{equation} \label{eq:temp1.ET}
\xymatrix{o \ar[r] & K \ar[r]^-u & X \ar[r] & S \ar[r] & o}
\end{equation}
with $X \in \mathcal{X}$. Let $\lambda$ be a regular cardinal such that $K$ is generated by a subset of cardinality less than $\lambda$.
By  transfinite induction on all ordinals $\alpha \leq \lambda$, we construct a 
$\frak B$-proper continuous chain $\left \{ P_\alpha \right \}_{\alpha \leq \lambda}$ of $\Lambda$-modules such that $P_{\alpha+1}/P_\alpha$ is isomorphic to a 
direct sum of copies of $S$ for every $\alpha < \lambda$. Let $P_0 = M$ and assume that $P_\alpha$ has already been constructed for an ordinal
$\alpha < \lambda$. Let $I_\alpha : = \hom_\Lambda (K,P_\alpha)$ and consider the pushout diagram
\begin{equation} \label{eq:temp21}
\begin{array}{c}
\xymatrix{
o \ar[r] & K^{(I_\alpha)} \ar[r]^-{u^{(I_\alpha)}} \ar[d]_-{\nu^K_{P_\alpha}} & X^{(I_\alpha)} \ar[d]^-{\varphi_\alpha} \ar[r] & S^{(I_\alpha)} \ar[r] \ar@{=}[d] & o \\
o \ar[r] & P_\alpha \ar[r]^-{i_\alpha} & D \ar[r] & S^{(I_\alpha)} \ar[r] & o
}
\end{array}
\end{equation}
wherein $\nu^K_{P_\alpha} : K^{(I_\alpha)} \arrow P_{\alpha}$ is the natural evaluation map and $u^{(I_\alpha)}$ is the map naturally induced by $u$. 
Note that since the upper row of the diagram is $\frak B$-exact, the lower row is also $\frak B$-exact by Observation~\ref{pullback.of.Bexact}.
Let $P_{\alpha+1}:=D$ and if, furthermore, $\alpha \leq \lambda$ is a limit ordinal and $P_\beta$ is already constructed for every $\beta < \alpha$, then let $P_\alpha = \bigcup_{\beta < \alpha} P_\beta$. The family $\{ P_\alpha \}_{\alpha \leq \lambda}$ is the desired $\frak B$-proper continuous chain.\\

\noindent \textsc{Step~2.} Let $P=\bigcup_{\alpha \leq \lambda} P_\alpha$, and for every $\alpha < \lambda$ 
let $\xymatrix{\eta_\alpha : P_\alpha \ar@{^{(}->}[r] & P}$ and $\xymatrix{i_\alpha : P_\alpha \ar@{^{(}->}[r] & P_{\alpha+1}}$ be the inclusions.
We claim that the short exact sequence
\begin{equation} \label{eq:M->P}
\xymatrix{ o \ar[r] & M \ar[r]^-{\eta_0} & P \ar[r]^-p & N \ar[r] & o } 
\end{equation}
is $\mathfrak{B}$-exact. In order to prove this claim, it suffices to show that every $\Lambda$-homomorphism $f : M \arrow Y$ with $Y \in \mathcal{Y}$ 
factors through $\eta_0$ via some $\Lambda$-homomorphism $g : P \arrow Y$, and finding such a $g$ amounts to constructing a family
$\{P_\alpha \stackrel{g_\alpha}{\arrow} Y \}_{\alpha \leq }$ of $\Lambda$-homomorphisms with the following \emph{compatibility properties}:
\begin{enumerate} \renewcommand{\theenumi}{\alph{enumi}}
\item $g_{\alpha+1} \mid_{P_\alpha} = g_\alpha$ for all $\alpha \leq \lambda$;

\item $g_\alpha \mid_M = f$ for all $\alpha \leq \lambda$.
\end{enumerate}
Let $g_0 := f$ and assume that $g_\alpha$ is already constructed for some $\alpha < \lambda$ . Since $Y \in \mathcal{Y}$ and 
 the inclusion $\xymatrix{i_\alpha : P_\alpha \ar@{^{(}->}[r] & P_{\alpha+1}}$ is a $\mathfrak{B}$-monomorphism in the sense of~\ref{df:B-exactness}
by construction in \textsc{Step~1}, there exists a $\Lambda$-homomorphism $g_{\alpha+1} : P_{\alpha+1} \arrow Y$ making the diagram
\[\xymatrix{ 
P_\alpha \ar[d]_-{g_\alpha} \ar@{^{(}->}[r] & P_{\alpha+1} \ar@{-->}[dl]^-{g_{\alpha+1}} \\
Y &
 }\]
commute. If, furthermore, $\alpha \leq \lambda$ is a limit ordinal and the family $\{P_\beta \stackrel{g_\beta}{\arrow} Y \}_{\beta < \alpha}$ is already
constructed with the compatibility properties (a) and (b) mentioned above, then we can simply let $g_\alpha = \bigcup_{\beta < \alpha} g_\beta$.
Now the family $\{ P_{\alpha} \stackrel{g_\alpha}{\arrow} Y \}_{\alpha \leq \lambda}$ has the desired properties (a) and (b), and 
the $\Lambda$-homomorphism $g=\bigcup_{\alpha \leq \lambda} g_\alpha$ satisfies $g \circ \eta_0 = f$. Thus, the sequence~\eqref{eq:M->P} is
$\mathfrak{B}$-exact.\\

\noindent \textsc{Step~3.} We finish the proof of part~(i) by showing that in the $\mathfrak{B}$-exact sequence~\eqref{eq:M->P}
we have $P \in \mathcal{S}^\bperp$ and $N$ is $\mathfrak{B}$-properly $\mathcal{S}$-filtered.
Applying $\hom_\Lambda (-,P)$ to~\eqref{eq:temp1.ET} yields the exact sequence
\[\xymatrix{ \hom_\Lambda (X,P) \ar[r]^-{u^\ast} & \hom_\Lambda (K,P) \ar[r] & \Ext^1_\frak B (S,P) \ar[r] & o }\: .  \]
Therefore, $P \in S^\bperp$ amounts to proving that every $\Lambda$-homomorphism $f : K \arrow P$ factors 
through $u$. In order to prove this, notice first that since $K$ is generated by a subset of cardinality less than $\lambda$
and $\lambda$ is regular, $f(K) \subseteq P_\alpha$ for some $\alpha < \lambda$. Let $f_\alpha : K \arrow P_\alpha$
be the $\Lambda$-homomorphism obtained by restricting the codomain of $f$ to $P_\alpha$. Thus, $f=\eta_\alpha \circ f_\alpha$ and note that
 $f_\alpha$ factors through $\nu^K_{P_\alpha} : K^{(I_\alpha)} \arrow P_\alpha$ 
via the natural injection $\imath_{f_\alpha} : K \arrow K^{(I_\alpha)}$ corresponding to $f_\alpha \in I_\alpha$. 
Now the commutative diagram
\[\xymatrix{
K \ar[r]^-{\imath_{f_\alpha}} \ar[d]_-{u} & K^{(I_\alpha)} \ar[r]^-{\nu^K_{P_\alpha}} \ar[d]_-{u^{(I_\alpha)}} & 
P_\alpha \ar[r]^-{\eta_\alpha} \ar[d]_-{i_\alpha} & P \ar@{=}[d] \\
X \ar[r]_-{j_{f_\alpha}} & X^{(I_\alpha)} \ar[r]_-{\varphi_\alpha} & P_{\alpha+1} \ar[r]_-{\eta_{\alpha+1}} & P
}\]
wherein $j_{f_\alpha} : X \arrow X^{(I_\alpha)}$ is the canonical injection corresponding to $f_\alpha \in I_\alpha$,
shows that $f$ factors through $u$ via $\eta_{\alpha+1} \circ \varphi_{\alpha} \circ j_{f_\alpha}$.  
Consequently, $P \in S^{\bperp}$.
Finally, for every $\alpha \leq \lambda$ let $N_\alpha : = P_\alpha / M$ 
and note that for any ordinal $\alpha < \lambda$, $N_{\alpha+1} / N_{\alpha} \iso P_{\alpha+1} / P_{\alpha} \iso S^{(I_\alpha)} \in {}^\bperp (S^\bperp)$.
Furthermore, if $\pi_\alpha : P_\alpha \arrow N_\alpha$ is the natural surjection and 
$\bar{i_\alpha} : N_\alpha \arrow N_{\alpha+1}$ is the injection induced by $i_\alpha : P_{\alpha} \arrow P_{\alpha+1}$
for every ordinal $\alpha < \lambda$, then the pushout diagram
\[\xymatrix{
o \ar[r] & P_{\alpha} \ar[d]_-{\pi_\alpha} \ar[r]^-{i_\alpha} & P_{\alpha+1} \ar[r] \ar[d]_-{\pi_{\alpha+1}} &
 P_{\alpha+1} / P_{\alpha} \ar[r] \ar[d]^-{\cong} & o \\
o \ar[r] & N_{\alpha} \ar[r]^-{\bar{i_\alpha}}  & N_{\alpha+1} \ar[r] & N_{\alpha+1} / N_{\alpha} \ar[r] & o
}\]
by virtue of~\ref{pullback.of.Bexact} implies that the short exact sequence
\[\xymatrix{ o \ar[r] & N_\alpha \ar[r]^-{\bar{i_{\alpha}}} & N_{\alpha+1} \ar[r] & N_{\alpha+1} / N_{\alpha} \ar[r] & o }\]
is $\frak B$-exact for every $\alpha < \lambda$. Consequently, $N$ is $\frak B$-properly $\mathcal{S}$-filtered by $\{N_\alpha\}_{\alpha \leq \lambda}$.
Thus, the sequence \eqref{eq:M->P} is the desired $\mathfrak{B}$-exact sequence in the assertion and the proof of part~(i) of the theorem is complete.\\

Part~(ii): Completeness of $\big ( {}^\bperp (\mathcal{S}^\bperp) , \mathcal{S}^\bperp \big )$ follows immediately from part~(i).
As for the second assertion, let $\mathcal{Z}$ be the class of all direct summands of direct sums of $\mathcal{S}^+$-filtered modules.
Notice first that $\mathcal{Z} \subseteq {}^\bperp (\mathcal{S}^\bperp)$ by Relative Eklof Lemma~\ref{Eklof.lemma.Rel}. On the other hand, 
for every $M \in {}^\bperp (\mathcal{S}^\bperp)$ consider, by admissibility of $\mathfrak{B}$, an exact a short $\frak B$-exact sequence
\[\xymatrix{ o \ar[r] & K \ar[r]^-u & X \ar[r] & M \ar[r] & o }\]
where $X \in \mathcal{X}$ and choose, by part~(i), an exact short $\frak B$-exact sequence
\[\xymatrix{ o \ar[r] & K \ar[r]^-v & P \ar[r] & N \ar[r] & o }\]
where $P \in \mathcal{S}^\bperp$ and $N$ is $\frak B$-properly $\mathcal{S}$-filtered. Now form the pushout of $u$ and $v$ to obtain 
the following commutative diagram with exact rows and columns:
\begin{equation} \label{eq:ET221}
\begin{array}{c}
\xymatrix{
         & o \ar[d]            & o \ar[d] &                            &  \\
o \ar[r] & K \ar[r]^-u \ar[d]_-{v}       & X \ar[r] \ar[d] & M \ar[r] \ar@{=}[d] & o \\
o \ar[r] & P \ar[r]  \ar[d]          & Z \ar[r] \ar[d] & M \ar[r] & o \\
         &  N \ar@{=}[r] \ar[d] & N      \ar[d] &           & \\ 
         & o                    & o       &           &
}\end{array}
\end{equation}
Note that the lower row and the right-hand side column of the diagram are $\frak B$-exact by~\ref{pullback.of.Bexact}.
Since $M \in {}^\bperp (\mathcal{S}^\bperp)$ and $P \in \mathcal{S}^\bperp$, the middle row of the diagram splits so that
$M$ is a direct summand of $Z$. On the other hand, since $N$ is $\frak B$-property $\mathcal{S}$-filtered, the middle column of the
diagram shows that $Z \in \mathcal{Z}$. Therefore, $M \in \mathcal{Z}$ and so $\mathcal{Z} = {}^\bperp (\mathcal{S}^\bperp)$.
\eproof

\brem \label{ET.rem}
Note that in~\eqref{eq:ET221}, $P \in {}^{\bperp} (\mathcal{S}^{\bperp}) \cap \mathcal{S}^{\bperp}$. Therefore, we may refine part~(ii) of Theorem~\ref{thm:eklof.trlifaj*} as follows: \textit{$M \in {}^{\bperp} (\mathcal{S}^{\bperp})$ if and only if there exists
$N \in {}^{\bperp} (\mathcal{S}^{\bperp}) \cap \mathcal{S}^{\bperp}$ such that $M \amalg N$ is $\mathcal{S}^+$-filtered}.
\erem

\bcoro \label{rel.cotor.T} 
\label{lem:T.class.countable*}
Let $T$ be a $\Lambda$-module and $X_\bullet$ be a left  $\mathcal{X}$-resolution of $T$. If $\mathcal{C}$ is the set of all
the relative syzygy modules of $X_\bullet$, then $T^{\bperp_\infty} = \mathcal{C}^{\bperp}$ and
$\mathfrak{C}_T : = \big ( {}^\bperp (T^{\bperp_\infty}) , T^{\bperp_\infty} \big )$ is a complete hereditary $\mathfrak{B}$-cotorsion pair 
cogenerated by $\mathcal{C}$.
\ecoro

\bproof
The equality $T^{\bperp_\infty} = \mathcal{C}^\bperp$ holds by the fact that for every $n \geq 1$,
\[ \Ext^{n+1}_\frak B \big (T , - \big) \iso \Ext^1_\frak B \big ( \Omega_{n} (X_\bullet) , - \big )\]
and this immediately implies  $\mathfrak{C}_T$ is the cotorsion pair cogenerated by the set $\mathcal{C}$ and hence complete by Theorem~\ref{thm:eklof.trlifaj*}.
Finally, since the class $T^{\bperp_\infty}$ is $\frak B$-coresolving by~\ref{obs:*inf.classes}, one concludes by Proposition~\ref{prop:*rogas} that the
$\mathfrak{B}$-cotorsion pair $\mathfrak{C}_T$ is hereditary.
\eproof

\brem \label{rem:size.of.S}
A couple of remarks are in order regarding the size of the class $\mathcal{S}^+$ and the structure of ${}^{\bperp} (\mathcal{S}^{\bperp})$ in Theorem~\ref{thm:eklof.trlifaj*}:
\begin{enumerate}
\item If the class $\mathcal{X}$ has a \emph{large representation generator}, i.e. there exists a $\Lambda$-module $G$ with $\mathcal{X} = \Add_\Lambda (G)$, then we can let $\mathcal{S}^{+}:=\mathcal{S} \cup \{G\}$ in Theorem~\ref{thm:eklof.trlifaj*}: Indeed,
if $\mathcal{X} = \Add_\Lambda (G)$, then  for every $\Lambda$-module $M$ there exists an exact $\frak B$-exact sequence of the form
\[\xymatrix{ o \ar[r] & K \ar[r] & G^{(\kappa)} \ar[r] & M \ar[r] & o }\]
for some cardinal number $\kappa$, and an argument similar to the proof of part~(ii) of Theorem~\ref{thm:eklof.trlifaj*} shows that
$M$ is a direct summand of a $\mathfrak{B}$-properly $\mathcal{S}^+$-filtered module where $\mathcal{S}^{+}=\mathcal{S} \cup \{G\}$.

\item The situation described in part~(i) happens in the following important case: Let $\Lambda$ be an artin algebra and assume 
that $\mathcal{X} = \dirlim \mathcal{X}^\fin$ where $\mathcal{X}^\fin$ is of finite representation type.
Then by~\cite[Theorem~3.1]{beligiannis.CMtype}, if $\mathcal{X}^\fin$ is contravariantly finite and resolving in $\mod (\Lambda)$, then every module $\mathcal{X}$ is a direct sum of modules in $\mathcal{X}^\fin$. Note that in this case $\mathcal{X}$ is closed under direct sums (cf.~\cite[Lemma 2.13]{GT}) and hence 
if $G$ is a representation generator for $\mathcal{X}^\fin$, then $\mathcal{X} = \Add_\Lambda (G)$.
\end{enumerate}
\erem

\bthm \label{thm:size.of.S} 
If $\Lambda$ is a CM-finite virtually Gorenstein artin algebra with $\Gp_0 = \add_\Lambda (H)$,
then for any set $\mathcal{S}$ of $\Lambda$-modules the class ${}^{\Gperp} (\mathcal{S}^{\Gperp})$ consists precisely of direct 
summands of $G$-properly $\mathcal{S}^+$-filtered $\Lambda$-modules where $\mathcal{S}^+ = \mathcal{S} \cup \{ H \}$.
\ethm

\bproof
Since $\Lambda$ is virtually Gorenstein, it follows from~\cite[Theorem~8.2]{beligiannis.CM} that $\GP_0 = \dirlim \Gp_0$. 
The class $\Gp_0$ is of finite representation type, contravariantly finite, and resolving in $\mod (\Lambda)$. 
Thus the result follows from Remark~\ref{rem:size.of.S}.
\eproof

\section{Gorenstein Tilting Modules and Their Basic Properties}

In this section, we apply the results of the previous sections to study ``infinitely Gorenstein tilting modules'' and their ``$G$-approximations''
via  $G$-cotorsion pairs over virtually Gorenstein rings. We start with the definition of an ``infinitely generated Gorenstein tilting module'' which is modeled on the definition of a infinitely generated tilting module~\cite{hugel.coel} and Auslander-Solberg's definition~\cite[Section~3]{ASII} of a finitely generated relative tilting module; cf. also~\cite[Definition~3.2]{Gtilting}.

\bdes{Gorenstein Tilting Modules} \label{df:G.tilting} 
Let $\Lambda$ be a virtually Gorenstein ring  and $n \geq 0$ be an integer. A $\Lambda$-module $T$ is called a \emph{Gorenstein $n$-tilting module}
if the following properties are satisfied:
\begin{description}
\item[(G.T1)] $\Gpd_\Lambda (T) \leq n$;

\item[(G.T2)] For any cardinal number $\kappa$, $\Gext^{\geqslant 1}_\Lambda (T,T^{(\kappa)})=o$;

\item[(G.T3)] For all $G \in \GP_0$, there exists a $G$-exact sequence
\[\xymatrix{ o \ar[r] & G \ar[r] & T_0 \ar[r] & \cdots \ar[r] & T_m \ar[r] & o }\]
where $T_i \in \Add_\Lambda (T)$ for all $0 \leq i \leq m$.
Here $\Add_\Lambda (T)$ denotes the class of all $\Lambda$-modules which are direct summands of direct sums of copies of $T$.
\end{description}
A $\Lambda$-module satisfying conditions (G.T1)--(G.T2) above is called  \emph{Gorenstein $n$-exceptional},
and a Gorenstein $n$-exceptional $\Lambda$-module $T$ is called \emph{partial Gorenstein $n$-tilting} if $T^{\Gperp_\infty}$ is closed under direct sums.
\edes

The next three fundamental lemmas, which are the Gorenstein versions of Lemma 2.2--2.4 in~\cite{hugel.coel},
are of particular importance in investigation of Gorenstein tilting modules.

\blem[cf.~{\cite[Lemma~2.2]{hugel.coel}}] \label{lem:tilting.class}
Let $\Lambda$ be a virtually Gorenstein ring and $T$ be a $\Lambda$-module. If $n:=\Gpd_\Lambda (T)<+\infty$, then any $\Lambda$-module has a finite $G$-exact $T^{\Gperp_\infty}$-coresolution of length $n$, and ${}^{\Gperp} (T^{\Gperp_\infty}) \subseteq \mathcal{GP}_n$.
\elem

\bproof
For any $\Lambda$-module $M$, take an $n$-step right $\GI_0$-resolution
\[\xymatrix{ \delta := o  \ar[r] & M \ar[r] & I_0 \ar[r] & \cdots \ar[r] & I_{n-1} \ar[r] & Q_n \ar[r] & o } \: . \]
Since $\Gpd_\Lambda (T)=n$, it follows from dimension shifting that $Q_n \in T^{\Gperp_\infty}$. Therefore, $\delta$ is
a $G$-exact $T^{\Gperp_\infty}$-coresolution of $M$, and consequently every $\Lambda$-module has a finite $G$-exact $T^{\Gperp_\infty}$-coresolution
of length $n$. Assume now that $A \in {}^{\Gperp} (T^{\Gperp_\infty})$.
For any $\Lambda$-module $M$ take a $G$-exact $T^{\Gperp_\infty}$-coresolution of length $n$, as $\delta$ above, and note that
$\Gext^{n+1}_\Lambda (A,M) \iso \Gext^1_\Lambda (A,Q_n)$ by Dimension Shifting~\ref{lem:Dim.shifting}.
But $\Gext^1_\Lambda (A,Q_n)=o$ because $Q_n \in T^{\Gperp_\infty}$.
Therefore, $\Gpd_\Lambda (A) \leq n$ by~\ref{obs:BExt}. Consequently, the inclusion
${}^{\Gperp} (T^{\Gperp_\infty}) \subseteq \mathcal{GP}_n$ holds.
\eproof

Before we state the next lemma, we need to fix a notation and define the notions of ``finendo'' and ``endo-finite'' modules.

\bnotation \label{note:Gen.GP}
Given a module $T$ over a  virtually Gorenstein ring $\Lambda$, we denote by $\Gen_\GP (T)$ the class of all $\Lambda$-modules $M$ such that
there exists a $G$-epimorphism $T^{(\kappa)} \arrow M$ for some cardinal $\kappa$.
\enotation

\bdf \label{finendo} 
\label{endofinite}
Recall that a $\Lambda$-module $M$ is said be \emph{finendo} if it is  finitely generated as a module over $\End_\Lambda (T)^\op$,
and it is said to be \emph{endo-finite} if it is of finite length as a module over $\End_\Lambda (T)^\op$.
Clearly, every finitely generated module over an artin algebra is endo-finite.
\edf

By the following lemma every Gorenstein tilting module over a virtually Gorenstein ring is finendo. The lemma is indeed the Gorenstein analogue of~\cite[Lemma~2.3]{hugel.coel}

\blem \label{lem:tilting->Gen} 
Let $\Lambda$ be a virtually Gorenstein ring and $T$ be a $\Lambda$-module.
If $T$ satisfies (G.T2) and (G.T3), then $T^{\Gperp_\infty} \subseteq \Gen_{\mathcal{GP}} (T)$ and $T$ is finendo.
\elem

\bproof
Let $M \in T^{\Gperp_\infty}$ and consider an exact sequence
\[\xymatrix{ o \ar[r] & K \ar[r] & G \ar[r]^-f & M \ar[r] & o }\]
where $f$ is a $\GP_0$-precover of $M$. On the other hand, there exists by (G.T3) a $G$-exact sequence
\[\xymatrix{
o \ar[r] & G \ar[r]^-g & T_0 \ar[r] & \cdots \ar[r] &
T_m \ar[r] & o
}\]
where each $T_i$ belongs to $\Add_\Lambda (T)$ for all $0 \leq i \leq m$. Therefore, we can form the following push-out diagram
\[\xymatrix{
     &          & o \ar[d]            & o \ar[d] & \\
o \ar[r] & K \ar[r] \ar@{=}[d] & G \ar[r]^-f \ar[d]_-g & M \ar[r] \ar[d] & o \\
o \ar[r] & K \ar[r] & T_0 \ar[r] \ar[d] & 
D \ar[r] \ar[d] & o \\
&  & \coker (g) \ar@{=}[r] \ar[d] & \coker (g) \ar[d] & \\
&   &  o  & o  &
}\]
Since the left-hand side column of the diagram is $G$-exact and
$G$-exact sequences are closed under push-out by~\ref{pullback.of.Bexact}, the right-hand side column is also $G$-exact. It now follows, say by 
the Snake Lemma, that the lower row of the diagram is $G$-exact and so $D \in \Gen_{\mathcal{GP}} (T)$.
It then follows from dimension shifting~\ref{lem:Dim.shifting} that $\coker (g) \in {}^{\Gperp} (T^{\Gperp_\infty})$ 
and so the right-hand side column of the diagram splits. Consequently, $M \in \Gen_{\mathcal{GP}} (T)$ and this finishes the proof of the inclusion $T^{\Gperp_\infty} \subseteq \Gen_{\mathcal{GP}} (T)$.

In order to prove that $T$ is finendo, note that there exists by (G.T3) a $G$-exact sequence
\[\xymatrix{
o \ar[r] & \Lambda \ar[r]^-f & T_0 \ar[r] & \cdots \ar[r] & T_m \ar[r] & o \: .
}\]
It is easily seen by (G.T2) and Dimension Shifting~\ref{lem:Dim.shifting} that the above sequence remains exact under
$\hom_\Lambda \big ( - , \Add_\Lambda (T) \big )$.
In particular, the map $f : \Lambda \arrow T_0$ is an $\Add_\Lambda (T)$-preenvelope of $\Lambda$ which amounts to $T$ being finendo 
by~\cite[Proposition~1.2]{tilting.preenv}.
\eproof

\blem[cf.~{\cite[Lemma~2.4]{hugel.coel}}] \label{lem:AddT}
Let $\Lambda$ be a virtually Gorenstein ring and $T$ be a $\Lambda$-module.
Assume that $T$ satisfies (G.T2) in~\ref{df:G.tilting} and $T^{\Gperp_\infty} \subseteq \Gen_{\mathcal{GP}} (T)$. Then:
\begin{enumerate}
\item If $f : M \arrow X$ is an $\Add_\Lambda (T)$-precover of $X \in T^{\Gperp_\infty}$, then the sequence
\[\xymatrix{ o \ar[r] & \ker (f) \ar[r] & M \ar[r]^-f & X \ar[r] & o  }\]
is $G$-exact with $\ker (f) \in T^{\Gperp_\infty}$. Consequently, $\Add_\Lambda (T) = T^{\Gperp_\infty} \cap {}^{\Gperp} (T^{\Gperp_\infty})$.

\item For every $X \in {}^{\Gperp} (T^{\Gperp_\infty})$, there exists a $G$-exact sequence
\[\xymatrix{o \ar[r] & X \ar[r]^-f & M \ar[r] & X' \ar[r] & o}\] where $M \in \Add_\Lambda (T)$ and $X' \in {}^{\Gperp} (T^{\Gperp_\infty})$.
In particular, $f : X \arrow M$ is an $\Add_\Lambda (T)$-preenvelope of $X$.

\item Every map $\xymatrix{A \ar[r]^-f &X}$ with $A \in {}^{\Gperp} (T^{\Gperp_\infty})$ and $X \in T^{\Gperp_\infty}$ factors through some module in $\Add_\Lambda (T)$. 
\end{enumerate}
\elem

\bproof
Part~(i): Since by the hypothesis $T^{\Gperp_\infty} \subseteq \Gen_\GP (T)$, there exists a short $G$-exact sequence of the form
\[\xymatrix{\delta :=o \ar[r] & K \ar[r] & T^{(\kappa)} \ar[r]^-g & X \ar[r] & o } \: . \]
for some cardinal $\kappa$. Since $g$ is surjective and $g$ factors through $f$, $f$ must be sujective and 
we can form the following push-out diagram.
\[\begin{array}{c}
\xymatrix{
\delta = o \ar[r] & K \ar[r] \ar[d] & T^{(\kappa)} \ar[r]^-g \ar[d] & X \ar[r] \ar@{=}[d] & o \\
\eta = o \ar[r] & \ker (f) \ar[r] & M \ar[r]^-f & 
X \ar[r] & o
}\end{array}  \]
Since $\delta$ is $G$-exact, it now follows from~\ref{pullback.of.Bexact} that the sequence $\eta$ is $G$-exact. On the other hand, since $f$ is 
an $\Add_\Lambda (T)$-precover of $X$, the sequence $\eta$ is  $\hom_\Lambda (T,-)$-exact and this in conjunction with $M,X \in T^{\Gperp_\infty}$ implies $\ker (f) \in T^{\Gperp_\infty}$. Finally, note that  if $X \in {}^{\Gperp} (T^{\Gperp_\infty})$, the sequence $\eta$ splits and hence $X \in \Add_\Lambda (T)$.
Consequently, ${}^{\Gperp} (T^{\Gperp_\infty}) \cap T^{\Gperp_\infty} \subseteq \Add_\Lambda (T)$. The reverse inclusion also holds by the fact that
$T$ satisfies (G.T2). Thus, $\Add_\Lambda (T) = {}^{\Gperp} (T^{\Gperp_\infty}) \cap T^{\Gperp_\infty}$.

Part~(ii): By Theorem~\ref{thm:eklof.trlifaj*} there exists a $G$-special $T^{\Gperp_\infty}$-preenvelope
$f : X \arrow M$ of $X$. Let $X' : = \coker (f)$ and consider the short $G$-exact sequence
\[\xymatrix{\delta := o \ar[r] & X \ar[r]^f & M \ar[r] & X' \ar[r] & o} \:  . \]
Since $X,X' \in {}^{\Gperp} (T^{\Gperp_\infty})$ and the class ${}^{\Gperp} (T^{\Gperp_\infty})$ is closed under $G$-extensions by Observation~\ref{obs:*perp},
we have $M \in {}^{\Gperp} (T^{\Gperp_\infty}) \cap T^{\Gperp_\infty}=\Add_\Lambda (T)$ by part~(i) of the lemma. Since $\delta$ is 
$G$-exact and $\Gext^1_\Lambda (X',T)=o$, the sequence remains exact under $\hom_\Lambda \big ( - , \Add_\Lambda (T) \big )$ and hence
$f$ is an $\Add_\Lambda (T)$-preenvelope of $X$.

Part~(iii): There exists by part~(i) a $G$-exact sequence of the form
\[\xymatrix{ o \ar[r] & X' \ar[r] & M \ar[r]^-g & X \ar[r] & o } \: . \]
where $g$ is an $\Add_\Lambda (T)$-precover of $X$. Since $A \in {}^{\Gperp} (T^{\Gperp_\infty})$ and $X' \in T^{\Gperp_\infty}$,
we have $\Gext^1_\Lambda (A,X')=o$ and so the sequence
\[\xymatrix{
o \ar[r] & \hom_\Lambda (A, X') \ar[r] & 
\hom_\Lambda (A , M ) \ar[r]^-g & \hom_\Lambda (A, X) \ar[r] & o }\]
is exact. Consequently, $f$ factors through $g$.
\eproof

Putting lemmas~\ref{lem:tilting.class}--\ref{lem:AddT} together we obtain the following useful characterization of Gorenstein tilting modules.

\bprop \label{prop:Gtilting.char}
Let $T$ be a Gorenstein $n$-exceptional module over a virtually Gorenstein ring $\Lambda$.
The following statements are equivalent:
\begin{enumerate}
\item  $T$ is Gorenstein $n$-tilting;

\item $T^{\Gperp_\infty} \subseteq \Gen_{\mathcal{GP}} (T)$;

\item For every $X \in {}^{\Gperp} (T^{\Gperp_\infty})$ there exists a $G$-exact sequence
\[\xymatrix{o \ar[r] & X \ar[r]^-f & M \ar[r] & X' \ar[r] & o}\] where $M \in \Add_\Lambda (T)$ and $X' \in {}^{\Gperp} (T^{\Gperp_\infty})$;

\item Every element of ${}^{\Gperp} (T^{\Gperp_\infty})$ admits a $G$-exact $\Add_\Lambda (T)$-coresolution of finite length.
\end{enumerate}
\eprop

\bproof
The implication (i)$\implies$(ii) follows from Lemma~\ref{lem:tilting->Gen} and the implication (ii)$\implies$(iii) follows from Lemma~\ref{lem:AddT}.
As for the implication (iii)$\implies$(iv), let $A \in {}^{\Gperp} (T^{\Gperp_\infty})$ and by repetitive use of~(iii) construct a $G$-exact sequence
of the form
\[\xymatrix{ o \ar[r] & A \ar[r]^-{d_0} & T_0 \ar[r] & \cdots \ar[r]^-{d_{n}} & T_n \ar[r] & C_n \ar[r] & o }\]
such that for all $0 \leq i \leq n$, $T_i \in \Add_\Lambda (T)$ and $C_i : = \coker (d_{i})$ belongs to ${}^{\Gperp} (T^{\Gperp_\infty})$.
Let $C_{-1} : = A$ and note that $\Gpd_\Lambda (A) \leq n$ by  Lemma~\ref{lem:tilting.class}.
Thus, Dimension Shifting~\ref{lem:Dim.shifting} yields
\[ \Gext^1_\Lambda (C_n , C_{n-1})  \cong \Gext^{n+1}_\Lambda (A, C_{n-1})  = o \: . \]
Consequently, the short $G$-exact sequence
\[\xymatrix{ o \ar[r] & C_{n-1} \ar[r] & T_n \ar[r] & C_n \ar[r] & o }\]
splits and so $C_{n-1} \in \Add_\Lambda (T)$. Therefore, $A$ admits a $G$-exact $\Add_\Lambda (T)$-coresolution of length $n$.
Finally, if (iv) holds then $T$ satisfies property (G.T3) in~\ref{df:G.tilting} because $\GP_0 \subseteq {}^{\Gperp} (T^{\Gperp_\infty})$.
Thus, $T$ is a Gorenstein $n$-tilting module and this completes the proof.
\eproof

We are now going to discuss ``Gorenstein tilting classes'' which are the components of the $G$-cotorsion pair induced by a Gorenstein tilting module.

\bdes{Gorenstein Tilting Classes} \label{df:tilting.classes} \label{df:Gtilting.cotorsion.pair}
Given a Gorenstein $n$-tilting module $T$ over a virtually Gorenstein ring $\Lambda$,
let $\mathcal{L}_T : = {}^{\Gperp} (T^{\Gperp_\infty})$ and $\mathcal{R}_T : = T^{\Gperp_\infty}$.
The class $\mathcal{L}_T$ is called the \emph{left Gorenstein ($n$-)tilting class} of $T$, and similarly the class $\mathcal{R}_T$ is called 
the \emph{right Gorenstein ($n$-)tilting class} of $T$. The two classes form the complete hereditary $G$-cotorsion pair $(\mathcal{L}_T , \mathcal{R}_T)$ by~\ref{rel.cotor.T}, which is referred to as the \emph{(Gorensteion $n$-tilting) $G$-cotorsion pair induced by $T$}.
\edes

In the following proposition we characterize the elements of Gorenstein tilting classes associated with a Gorenstein tilting module $T$ in terms of their $\Add_\Lambda (T)$-(co)resolutions. This characterization should be regarded as an infinitely generated version of~\cite[Theorem~3.2]{ASII} in the Gorenstein setting and also the
``Gorenstein analogue'' of~\cite[Proposition 13.13]{GT}.

\bprop \label{prop:Gtilting.class}
Let $T$ be a Gorenstein $n$-tilting module  over a virtually Gorenstein ring $\Lambda$. If $(\mathcal{L}_T , \mathcal{R}_T)$
is the $G$-cotorsion pair induced by $T$, then:
\begin{enumerate}
\item The class $\mathcal{L}_T$ consists precisely of all $\Lambda$-modules with a $G$-exact $\Add_\Lambda (T)$-coresolution of length at most $n$.

\item The class $\mathcal{R}_T$ consists precisely of all $\Lambda$-modules which have a $G$-exact $\Add_\Lambda (T)$-resolution. 
In particular, $\mathcal{R}_T$ is closed under direct sums.

\item For every integer $m \geq 0$, every module in $\mathcal{R}_T \cap \mathcal{GP}_m$ has
a $G$-exact $\Add_\Lambda (T)$-resolution of length at most $m$. 
\end{enumerate}
\eprop

\bproof 
Part~(i): Since $\mathcal{L}_T$ is $G$-resolving by~\ref{obs:*inf.classes}, every module which admits a $G$-exact $\Add_\Lambda (T)$-coresolution of 
finite length belongs to $\mathcal{L}_T$. On the other hand, if $M \in \mathcal{L}_T$, then iterated $G$-special $\mathcal{R}_T$-preenvelopes of $M$ yield 
by part~(iii) of Proposition~\ref{prop:Gtilting.char} a $G$-exact sequence
\[\xymatrix{ \tau : = o \ar[r] & M \ar[r]^-{d_0} & T_0 \ar[r]^-{d_1} & \cdots \ar[r]^-{d_n} & T_n \ar[r]^-{d_{n+1}} & T_{n+1} \ar[r] & \cdots  } \]
where each $T_i$ belongs to $\Add_\Lambda (T)$. Let $C_{-1} : = M$, $C_i := \coker (d_i)$ for every $i \geq 0$,
and note that since $\tau$ is obtained from iterated $G$-special $\mathcal{R}_T$-preenvelopes, $C_i \in \mathcal{L}_T \subseteq \mathcal{GP}_n$ for all
$i \geq 0$, where the inclusion $\mathcal{L}_T \subseteq \mathcal{GP}_n$ follows from Lemma~\ref{lem:tilting.class}. It now follows 
from part~(ii) of Observation~\ref{obs:BExt} and Dimension Shifting~\ref{lem:Dim.shifting} that
\[ \Gext^1_\Lambda (C_n , C_{n-1} ) \cong \Gext^{n+1}_\Lambda (C_n ,M) =o \: . \]
Thus, the canonical short $G$-exact sequence
\[\xymatrix{ o \ar[r] & C_{n-1} \ar[r] & T_n \ar[r] & C_{n} \ar[r] & o }\]
splits which implies that $C_{n-1} \in \Add_\Lambda (T)$. Consequently, $M$ admits a $G$-exact $\Add_\Lambda (T)$-coresolution of length at most $n$.

Part~(ii): If $M \in \mathcal{R}_T$, then iterated $G$-special $\mathcal{L}_T$-precovers of $M$ yield by 
part~(i) of Lemma~\ref{lem:AddT} a $G$-exact $\Add_\Lambda (T)$-resolution of $M$. Conversely, assume $M$ is a $\Lambda$-module which admits a $G$-exact $\Add_\Lambda (T)$-resolution, 
say
\[\xymatrix{ \tau : = \cdots \ar[r] & T_1 \ar[r]^-{d_1} & T_0 \ar[r]^-{d_0} & M \ar[r]  & o } \: . \]
Let $K_{-1} : = M$ and $K_i : = \ker (d_i)$ for all $i \geq 0$. Since $T_i \in \Add_\Lambda (T) \subseteq \mathcal{R}_T$, it follows from Dimension Shifting~\ref{lem:Dim.shifting} that for every $i \geq 1$,
%
\[\xymatrix{ \Gext^i_\Lambda (T,M) \cong \Gext^{n+i}_\Lambda (T , K_{n-1} )  } \: .\]
But $\Gext^{n+i}_\Lambda (T , K_{n-1} ) =o$ because $\Gpd_\Lambda (T) \leq n$. 
Therefore, $M \in \mathcal{L}_T^{\Gperp} = \mathcal{R}_T$, and, consequently, $\mathcal{R}_T$ coincides with the class of all $\Lambda$-modules which admit a
$G$-exact $\Add_\Lambda (T)$-resolution. Since by~\ref{lem:dirlim.of.Gexact} the direct sum of $G$-exact $\Add_\Lambda (T)$-resolutions is again a $G$-exact $\Add_\Lambda (T)$-resolution, the class $\mathcal{R}_T$ is closed under direct sums. 

Part~(iii): 
For every $M \in \mathcal{R}_T \cap \GP_m$ use iterated $G$-special $\mathcal{L}_T$-precovers of $M$ and part~(i) of Lemma~\ref{lem:AddT} to 
construct a $G$-exact sequence of the form
\[\xymatrix{
o \ar[r] & K_{m} \ar[r] & T_m \ar[r]^-{d_m} & T_{m-1} \ar[r]^-{d_{m-1}} & \cdots \ar[r] & T_0 \ar[r]^-{d_0} & M \ar[r] & o } \]
where $T_i \in \Add_\Lambda (T)$ and $K_i : = \ker (d_i)$ for all $0 \leq i \leq m$. Let also $K_{-1} : = M$. 
Since $K_m \in \mathcal{R}_T$ and each $T_i$ belongs to $\mathcal{L}_T$, 
$\Gext^1_\Lambda (T_i , K_m)=o$ for all $0 \leq i \leq m$. Therefore, 
\[ \Gext^1_\Lambda (K_{m-1} , K_m) \cong \Gext^{m+1}_\Lambda (M,K_m)  \]
by Dimension Shifting~\ref{lem:Dim.shifting}.  But $\Gext^{m+1}_\Lambda (M,K_m)=o$ because $\Gpd_\Lambda (M) \leq m$. 
Consequently, the short $G$-exact sequence
\[\xymatrix{ o \ar[r] & K_m \ar[r] & T_m \ar[r] & K_{m-1} \ar[r] & o }\]
splits, which implies that $K_{m-1} \in \Add_\Lambda (T)$. Therefore, $M$ admits a $G$-exact $\Add_\Lambda (T)$-resolution of length $m$.
\eproof

\section{Gorenstein Tilting Modules over Virtually Gorenstein Algebras of Finite CM-type}

In this section we focus on Gorenstein tilting modules over virtually  Gorenstein CM-finite artin algebras. 
As we shall see in Proposition~\ref{prop:Gtilting.over.CMfinite}, the condition (G.T3) in Definition~\ref{df:G.tilting} 
can be relaxed to a somewhat weaker condition over CM-finite artin algebras, and this allows us to give a characterization of right 
Gorenstein tilting classes over these algebras analogues to H\"{u}gel-Coelho's characterization~\cite[Theorem~4.1]{hugel.coel} of right tilting 
classes in ``standard'' tilting theory; cf. Theorem~\ref{thm:Gtilting.class}. Some applications of this characterization will be discussed in
 sections~\ref{sec:applications}.

\blem \label{lem:H}
Let $\Lambda$ be a virtually Gorenstein ring such that $\GP_0 = \Add_\Lambda (H)$ for some $H \in \GP_0$.
%
\begin{enumerate}
\item Every $\Lambda$-module $M$ has a $\mathcal{GP}_0$-precover of the form $f : H^{(\alpha)} \arrow M$ for some cardinal number $\alpha$. 
If, furthermore, $M$ is finitely generated, then $\alpha$ can be taken finite.

\item Let $T$ be a $\Lambda$-module which satisfies (G.T2) in~\ref{df:G.tilting}. If $H$ has a $G$-exact $\Add_\Lambda (T)$-coresolution of finite length,
then $T^{\Gperp_\infty} \subseteq \Gen_{\mathcal{GP}} (T)$.
\end{enumerate}
\elem

\bproof
Part~(i): Let $f : G \arrow M$ be a $\GP_0$-precover of $M$. 
Since $G \in \Add_\Lambda (H)$, there exists a Gorenstein projective $\Lambda$-module $G'$ such that $G \amalg G' \cong H^{(\alpha)}$ for some cardinal $\alpha$.
It is now easy to see that the $\Lambda$-homomorphism $\tilde{f} : H^{(\alpha)} \arrow M$ induced by $f$ is a $\GP_0$-precover of $M$.
If, furthermore, $G$ is finitely generated, then $\Lambda$-homomorphism $\tilde{f}$ has a factorization of the form
$\xymatrix{H^{(\alpha)} \ar[r] & H^n \ar[r]^-{\tilde{f}|_{H^n}} & M}$ for some integer $n \geq 0$, and it then follows easily that
${\tilde{f}|_{H^n}} : H^n \arrow M$ is a $\GP_0$-precover of $M$.

Part~(ii): Let $M \in T^{\Gperp_\infty}$ and choose by part~(i) of the lemma a (necessarily surjective) $\GP_0$-precover of $M$
of the form $f : H^{(\alpha)} \arrow M$ for some cardinal $\alpha$. Let
\[\xymatrix{ o \ar[r] & H \ar[r]^-{d} & T_0 \ar[r] & \cdots \ar[r] & T_m \ar[r] & o }\]
be a $G$-exact $\Add_\Lambda (T)$-coresolution of $H$ of finite length. Let $C : = \coker (d)$ and denote 
by $d^{(\alpha)} : H^{(\alpha)} \arrow T^{(\alpha)}_0$ map induced by $d$. Note that $d^{(\alpha)}$ is a $G$-monomorphism 
by Lemma~\ref{lem:dirlim.of.Gexact}. Now form the pushout of $f$ and $d^{(\alpha)}$ to obtain the following commutative diagram with exact rows and
columns:
\[\xymatrix{
     &          & o \ar[d]            & o \ar[d] & \\
o \ar[r] & K \ar[r] \ar@{=}[d] & H^{(\alpha)} \ar[r]^-f \ar[d]_-{d^{(\alpha)}}  & M \ar[r] \ar[d] & o \\
o \ar[r] & K \ar[r] & T^{(\alpha)}_0 \ar[r] \ar[d] & 
D \ar[r] \ar[d] & o \\
&  & C^{(\alpha)} \ar@{=}[r] \ar[d] & C^{(\alpha)} \ar[d] & \\
&   &  o  & o  &
}\]
Since the left-hand side column of the diagram is $G$-exact, it follows form~\ref{pullback.of.Bexact} that the right-hand side column 
is also $G$-exact, and it then follows that the lower row of the diagram is $G$-exact. Thus, $D \in \Gen_{\mathcal{GP}} (T)$. 
On the other hand, since $M \in T^{\Gperp_\infty}$ and $C \in {}^{\Gperp} (T^{\Gperp})$, we have
$\Gext^1_\Lambda (C^{(\alpha)} , M)= o$. Thus, the right-hand side column of the diagram splits, and so $M \in \Gen_{\mathcal{GP}} (T)$.
\eproof

\brem \label{rem:CM.has.Add}
The condition ``$\GP_0 = \Add_\Lambda (H)$'' in Lemma~\ref{lem:H} is satisfied for a virtually Gorenstein CM-finite artin algebra. 
Indeed, if $\Lambda$ is a virtually Gorenstein CM-finite artin algebra and $H$ is a representation generator of $\Gp_0$, 
then $\GP_0 = \Add_\Lambda (H)$ by~\cite[Theorem 4.10]{beligiannis.CMtype}.
\erem

\bprop \label{prop:Gtilting.over.CMfinite}
Let $\Lambda$ be a virtually Gorenstein ring such that $\GP_0 = \Add_\Lambda (H)$ for some $H \in \GP_0$.
A $\Lambda$-module $T$ is Gorenstein $n$-tilting if and only if $T$ satisfies:
\begin{description}
\item[(G.T1)] $\Gpd_\Lambda (T) \leq n$;

\item[(G.T2)] $\Gext^{\geqslant 1}_\Lambda (T , T^{(\kappa)})=o$ for any cardinal $\kappa$;

\item[(G.T3$'$)] There exists a $G$-exact sequence of the form
\[\xymatrix{ o \ar[r] & H \ar[r] & T_0 \ar[r] & \cdots \ar[r] & T_m \ar[r] & o }\]
where $T_i \in \Add_\Lambda (T)$ for all $0 \leq i \leq m$.
\end{description}
\eprop

\bproof
If $T$ is a Gorenstein $n$-tilting module, then $T$ satisfies (G.T1), (G.T2) and (G.T3$'$) by the definition; cf. Definition~\ref{df:G.tilting}.
Conversely, assume that $T$ satisfies (G.T1), (G.T2) and (G.T3$'$). By Lemma~\ref{lem:H}, $T^{\Gperp_\infty} \subseteq \Gen_{\mathcal{GP}} (T)$
and hence $T$ is a Gorenstein $n$-tilting module by Proposition~\ref{prop:Gtilting.char}.
\eproof

We are now in a position to prove the main result of this section, which is a Gorenstein analogue of~\cite[Theorem~4.1]{hugel.coel}, characterizing 
Gorenstein tilting classes; cf. also~\cite[Theorem~3.24]{ASII}.

\bthm \label{thm:Gtilting.class}
Let $\Lambda$ be a virtually Gorenstein ring such that $\GP_0 = \Add_\Lambda (H)$ for some $H \in \GP_0$.
A class $\mathcal{C}$ of $\Lambda$-modules is a right Gorenstein $n$-tilting class
\begin{enumerate}
\item ${}^{\Gperp} \mathcal{C} \subseteq \mathcal{GP}_n$;

\item $\mathcal{C}$ is $G$-coresolving and it is closed under direct sums and summands;

\item $\mathcal{C}$ is $G$-special preenveloping.
\end{enumerate}
\ethm

\bproof
Assume first that $\mathcal{C}$ is a Gorenstein $n$-tilting class so that $\mathcal{C} = T^{\Gperp_\infty}$ for some Gorenstein $n$-tilting $\Lambda$-module $T$.
Then (i) follows from Lemma~\ref{lem:tilting.class}, (ii) follows from~\ref{obs:*inf.classes} and~\ref{prop:Gtilting.class},
and finally (iii) follows from Corollary~\ref{rel.cotor.T}.
Conversely, assume $\mathcal{C}$ satisfies the conditions (i)--(iii). The proof that $\mathcal{C}=T^{\Gperp_\infty}$ for some Gorenstein $n$-tilting $\Lambda$-module
$T$ proceeds in several steps:\\

\noindent \textsc{Step~1.} Let $A$ be a $\Lambda$-module. Since $\mathcal{C}$ is by the hypothesis $G$-special preenveloping,
there exists a $G$-exact sequence
\[\xymatrix{ o \ar[r] & A \ar[r]^-{d_0} & B_0 \ar[r] & \cdots \ar[r]^-{d_n} & B_n \ar[r] & C_n \ar[r] & o  }\]
where $B_i \in \mathcal{C}$ and $C_i := \coker (d_i) \in {}^{\Gperp} \mathcal{C} \subseteq \mathcal{GP}_n$ for all $0 \leq i \leq n$. 
Let $C_{-1}  := A$ and note that the inclusion ${}^{\Gperp} \mathcal{C} \subseteq \mathcal{GP}_n$ implies $\Gpd_\Lambda (C_n) \leq n$.
Thus, $\Gext^1_\Lambda (C_n , C_{n-1}) \cong \Gext^{n+1}_\Lambda (C_n , A)=o$  by Dimension Shifting~\ref{lem:Dim.shifting}.
Therefore, the short $G$-exact sequence
\[\xymatrix{ o \ar[r] & C_{n-1} \ar[r] & B_n \ar[r] & C_n \ar[r] & o }\]
splits, and so $C_{n-1} \in \mathcal{C}$ because $\mathcal{C}$ is closed under direct summands by the hypothesis.
Consequently, there exists a $G$-exact sequence
\[\xymatrix{
o \ar[r] & A \ar[r] & M_0 \ar[r] & \cdots \ar[r] & M_n \ar[r] & o
}\]
where $M_i \in \mathcal{C}$ and all the cokernels of the maps in the sequence belong to ${}^{\Gperp} \mathcal{C}$. 
In particular $M_n \in {}^{\Gperp} \mathcal{C} \cap \mathcal{C}$.\\

\noindent \textsc{Step~2.} By \textsc{Step~1}, there exists a $G$-exact sequence
\[\xymatrix{
\tau_\bullet := o \ar[r] & H \ar[r] & T_0 \ar[r] & \cdots \ar[r] & T_n \ar[r] & o
}\]
where $T_i \in \mathcal{C}$ and all the cokernels of the maps in the sequence belong to ${}^{\Gperp} \mathcal{C}$. Let $T : = \coprod_{i = 0}^n T_i$. We prove in the
subsequent steps that $T$ is a Gorenstein $n$-tilting module such that $\mathcal{C} = T^{\Gperp_\infty}$.\\

\noindent \textsc{Step~3.} $T$ is a Gorenstein $n$-tilting module: Since $H$ and the cokernels of maps in $\tau_\bullet$ all belong to ${}^{\Gperp} \mathcal{C}$,
it follows that $T_i \in \mathcal{C} \cap {}^{\Gperp} \mathcal{C}$ for all $0 \leq i \leq n$. In particular, 
each $T_i$ is of Gorenstein projective dimension at most $n$ by (i), which implies $\Gpd_\Lambda (T) \leq n$. Furthermore, since
$\mathcal{C}$ is closed under direct sums, $T^{(\kappa)} \in \mathcal{C} \cap {}^{\Gperp} \mathcal{C}$ for every cardinal $\kappa$,
which implies $\Gext^{\geqslant 1} (T , T^{(\kappa)}) = o$; note that ${}^{\Gperp} \mathcal{C}={}^{\Gperp_\infty
} \mathcal{C}$ by~\ref{obs:*inf.classes}. In view of the $G$-exact sequence $\tau_\bullet$, it now follows from Proposition~\ref{prop:Gtilting.over.CMfinite}
that $T$ is a Gorenstein $n$-tilting $\Lambda$-module.\\

\noindent \textsc{Step~4.} The equality $\mathcal{C} \cap {}^{\Gperp} \mathcal{C} = \Add_\Lambda (T)$ holds: 
Indeed, since $T \in \mathcal{C} \cap {}^{\Gperp} \mathcal{C}$, it follows from the closure properties of $\mathcal{C}$ and ${}^{\Gperp} \mathcal{C}$
that $\Add_\Lambda (T) \subseteq \mathcal{C} \cap {}^{\Gperp} \mathcal{C}$. As for the reverse inclusion, notice first that
$\mathcal{C} \cap {}^{\Gperp} \mathcal{C} \subseteq T^{\Gperp_\infty} \cap \mathcal{GP}_n$. Therefore, by Proposition~\ref{prop:Gtilting.class},
for every $M \in \mathcal{C} \cap {}^{\Gperp} \mathcal{C}$ there exists a $G$-exact sequence
\[\xymatrix{ o \ar[r] & U_n \ar[r] & \cdots \ar[r] & U_0 \ar[r]^-f & M \ar[r] & o }\]
where $U_i \in \Add_\Lambda (T)$ for all $0 \leq i \leq n$. Since $\Add_\Lambda (T) \subseteq \mathcal{C}$ and $\mathcal{C}$ is $G$-coresolving by the
hypothesis, we have $\ker (f) \in \mathcal{C}$, which implies that $f$ is an split epimorphism and so $M \in \Add_\Lambda (T)$. 
Consequently, $\mathcal{C} \cap {}^{\Gperp} \mathcal{C} = \Add_\Lambda (T)$. \\

\noindent \textsc{Step~5.} We finish the proof by showing that $\mathcal{C} = T^{\Gperp_\infty}$. Since $T \in {}^{\Gperp} \mathcal{C}$ and
$\mathcal{C}$ is $G$-coresolving, $\mathcal{C} \subseteq T^{\Gperp_\infty}$. As for the reverse inclusion, let $A \in T^{\Gperp_\infty}$ and
note that there exists by Step~1 a $G$-exact sequence
\[\xymatrix{ o \ar[r] & A \ar[r]^-{f_0} & M_0 \ar[r] & \cdots \ar[r]^-{f_n} & M_n \ar[r] & o  }\]
where $M_i \in \mathcal{C} \subseteq T^{\Gperp_\infty}$ for all $0 \leq i \leq n$,
and all the cokernels belong to ${}^{\Gperp} \mathcal{C}$. 
In particular, $M_n \in \mathcal{C} \cap {}^{\Gperp} \mathcal{C}=\Add_\Lambda (T) \subseteq {}^{\Gperp} (T^{\Gperp_\infty})$ by Step~4. Since $A \in T^{\Gperp_\infty}$ and $T^{\Gperp_\infty}$ is $G$-coresolving, all the cokernels also belong to $T^{\Gperp_\infty}$. In particular $\ker (f_n) \in T^{\Gperp_\infty}$.
Thus, $\Gext^1_\Lambda \big (M_n , \ker (f_n) \big ) =o$ and so $f_n$ is a split epimorphism. It then follows by induction that $f_0$ is a split monomorphism
and so $M \in \mathcal{C}$, because $\mathcal{C}$ is closed under direct summands by the hypothesis. Consequently, $\mathcal{C}=T^{\Gperp_\infty}$.
\eproof

The Gorenstein tilting module $T$ whose right Gorenstein tilting class is $\mathcal{C}$ in Theorem~\ref{thm:Gp.contra} is \textit{not} finitely generated in general,
and in fact, as we shall prove in the sequel, finitely generatedness of  $T$, under some mild conditions, characterizes contravariantly finiteness of
$({}^{\Gperp} \mathcal{C})^\fin$; cf. Corollary~\ref{coro:G.tilting.class}.
To this end, we first need the following lemma about the approximations provided by
finitely generated Gorenstein tilting modules.

\blem \label{lem:approx.fg.tilting}
Let $\Lambda$ be a CM-finite virtually Gorenstein artin algebra and $T$ be a finitely generated Gorenstein $n$-tilting module. For every finitely generated $\Lambda$-module $M$ there 
exist short $G$-exact sequences of the form
\begin{equation} \tag{$\dagger$}
\xymatrix{ o \ar[r] & M \ar[r] & U \ar[r] & X \ar[r] & o }
\end{equation}
and
\begin{equation} \tag{$\ddagger$}
\xymatrix{ o \ar[r] & Y \ar[r] & V \ar[r] & M \ar[r] & o }
\end{equation}
where $X$ and $V$ have finite $G$-exact $\add_\Lambda (T)$-coresolutions and $U,Y \in (T^{\Gperp_\infty})^\fin$.
In particular, the class $\mathcal{L}_T^\fin$ is contravariantly finite and $\mathcal{R}_T^\fin$ is covariantly finite. 
\elem

\bproof 
The proof proceeds in several steps:

\noindent \textsc{Step~1.} We first prove the assertion for the case $M \in (T^{\Gperp_\infty})^\fin$. Obviously, the short exact sequence
\[\xymatrix{ o \ar[r] & M \ar[r]^-{1_M} & M \ar[r] & o \ar[r] & o }\]
is of type ($\dagger$) in the statement of the lemma. As for existence of a short $G$-exact sequence of type ($\ddagger$) for $M$, notice that since $\Lambda$ is an artin algebra,
$\hom_\Lambda (T,M)$ is finitely generated as a module over the center $R$ of $\Lambda$. Let $\{f_1 , \ldots , f_n\}$ be a finite generating set of 
the finitely generated $R$-module $\hom_\Lambda (T,M)$ and let $f : T^n \arrow M$ be the map induced by $f_1 , \ldots , f_n$.
It is then easily seen that $f : T^n \arrow M$ is an $\Add_\Lambda (T)$-precover of $M$ and since $T^{\Gperp_\infty} \subseteq \Gen_\GP (T)$ by Lemma~\ref{lem:tilting->Gen},  the sequence
\[\xymatrix{ \delta:= o \ar[r] & K \ar[r] & T^n \ar[r]^-f & M \ar[r] & o }  \]
is $G$-exact. Now since $T^n , M \in T^{\Gperp_\infty}$, it follows from the long exact sequence of $\Gext$-modules
induced from the above short $G$-exact sequence that $\Gext^{\geqslant 2}_\Lambda (T,K)=o$.
Furthermore, since $f$ is an $\Add_\Lambda (T)$-precover, the sequence remains exact under $\hom_\Lambda (T,-)$. This in view of
$\Gext^1_\Lambda (T,T^n)=o$ implies $\Gext^1_\Lambda (T,K)=o$. Hence $K \in T^{\Gperp_\infty} \cap \mod (\Lambda)$, and so the sequence $\delta$
is of type ($\ddagger$) in the statement of the lemma.

\textsc{Step~2.} If $\xymatrixcolsep{0.5cm}\xymatrix{ o \ar[r] & M \ar[r] & U \ar[r]^-u & X \ar[r] & o }$ is a $G$-exact sequence of type ($\dagger$) in the statement of the lemma, then $M$ also admits a $G$-exact sequence of type ($\ddagger$): Indeed, the $\Lambda$-module $U \in (T^{\Gperp_\infty})^\fin$ 
sits by \textsc{Step~1} in a short $G$-exact sequence of the form
\[\xymatrix{ o \ar[r] & Y \ar[r] & T_0 \ar[r]^-{\tau} & U \ar[r] & o } \: , \]
where $T_0 \in \add_\Lambda (T)$ and $Y \in (T^{\Gperp_\infty})^\fin$.
Then we can form the commutative diagram
\[\xymatrix{
            & o   \ar[d]        & o    \ar[d]            &             &  \\
            & Y  \ar@{=}[r] \ar[d]           & Y \ar[d]                &             &  \\            
o \ar[r] & V \ar[d]   \ar[r] \ar@{}[dr]^-\hole|{(\ast)}        &  T_0 \ar[d]^-{\tau} \ar[r]^-{u \circ \tau}               &   X \ar@{=}[d] \ar[r]          & o  \\            
 o \ar[r] & M \ar[r] \ar[d] & U \ar[r]_-u \ar[d]  & X \ar[r]  & o \\
            & o           & o                &             & 
}\]
with exact rows and columns. Since the right-hand side column and the lower row of the diagram are $G$-exact,
it follows easily that the upper row of the diagram is $G$-exact. Furthermore, since the right-hand side column of the diagram is $G$-exact 
and the square~($\ast$) is pullback, it follows from~\ref{pullback.of.Bexact} that the left-hand side column is also $G$-exact.
Notice now  that $V$ has a $G$-exact $\add_\Lambda (T)$-coresolution of finite length because $X$ has such a coresolution and $T_0 \in \add_\Lambda (T)$.
Thus, the $G$-exact sequence $\xymatrixcolsep{0.5cm} \xymatrix{o \ar[r] & Y \ar[r] & V \ar[r] & M \ar[r] & o}$ from
the left-hand side column of the diagram above is a $G$-exact sequence of type ($\ddagger$) for $M$.

\textsc{Step~3.} Any finitely generated module $M$ admits a $G$-exact coresolution of finite length by
modules in $(T^{\Gperp_\infty})^\fin$: Indeed, since $\Lambda$ is a virtually Gorenstein artin algebra,
the class  $\GI_0^\fin$ is covariantly finite by~\cite[Theorem 8.2]{beligiannis.CM} and so we can take by a version of Wakamatsu Lemma in $\mod (\Lambda)$~\cite[Lemma~1.3]{AR.applications} an exact sequence
\[\xymatrix{ I_\bullet := o \ar[r] & M \ar[r] & I_0 \ar[r] & \cdots \ar[r] & I_n \ar[r] & C \ar[r] & o }\]
such that $I_0 , \ldots , I_n \in \GI^\fin_0$ and the cokernels of the maps in the exact sequence belong to ${}^{\perp} (\GI_0^\fin) \cap \mod (\Lambda)$.
But ${}^{\perp} (\GI_0^\fin) \cap \mod (\Lambda) \subseteq ({}^\perp \GI_0)^\fin = (\GP_0^\perp)^\fin$ by~\cite[Theorem~8.2-(x) and Proposition~8.13]{beligiannis.CM}. Consequently, $I_\bullet$ is $G$-exact. The modules $I_0 , \ldots , I_n$ obviously belong to $(T^{\Gperp_\infty})^\fin$, and it follows from $G$-exactness of $I_\bullet$ and dimension shifting~\ref{lem:Dim.shifting} that
$\Gext^{i}_\Lambda (T,C) \cong \Gext^{n+i+2}_\Lambda (T,M)=o$. Thus, $C \in (T^{\Gperp_\infty})^\fin$ and so
$I_\bullet$ is a coresolution of $M$ by $(T^{\Gperp_\infty})^\fin$.

\textsc{Step~4.} Thanks to \textsc{Step~3}, we can now complete the proof by induction on the length $m$ of a
$G$-exact $(T^{\Gperp_\infty})^\fin$-coresolution of $M$: For $m=0$ the assertion is already proved in \textsc{Step~1}.
Assume that $m>0$ and that the assertion holds for any finitely generated module with a $G$-exact
$(T^{\Gperp_\infty})^\fin$-coresolution of length less than $m$. Let 
\[\xymatrix{ o \ar[r] & M \ar[r]^-{d} & C_0 \ar[r] & \cdots \ar[r] & C_m \ar[r] & o  }\]
be a $G$-exact $(T^{\Gperp_\infty})^\fin$-coresolution of $M$ of length $m$. Let $N:=\coker (d)$ and consider the
short $G$-exact sequence $\xymatrix{o \ar[r] & M \ar[r]^-d & C_0 \ar[r]^-f & N \ar[r] & o}$. Since $N$ admits a
$G$-exact $(T^{\Gperp_\infty})^\fin$-coresolution of length $m-1$, it follows from the induction hypothesis that
$N$ sits in a short $G$-exact sequence of the form
\[\xymatrix{ o \ar[r] & K \ar[r] & X \ar[r]^-g & N \ar[r] & o }\]
where $K \in (T^{\Gperp_\infty})^\fin$ and $X$ has a finite $G$-exact $\add_\Lambda (T)$-coresolution.
Forming the pullback of $f$ and $g$ we obtain the commutative diagram
\[\xymatrix{
& & o \ar[d] & o \ar[d] &        \\                      
											 &       						  & K \ar@{=}[r] \ar[d] & K \ar[d]                & \\
 o \ar[r] & M \ar[r] \ar@{=}[d] & U \ar[r] \ar[d]         & X \ar[r] \ar[d]^-g   & o \\             
o \ar[r] & M \ar[r]                   & C_0 \ar[r]_-f \ar[d]     & N \ar[r]  \ar[d]  & o \\
& & o & o & 
}\]
with $G$-exact rows and columns by~\ref{pullback.of.Bexact}. Since $K , C_0 \in (T^{\Gperp_\infty})^\fin$ and the left-hand side column is $G$-exact,
we have $U \in (T^{\Gperp_\infty})^\fin$ and so the upper row of the diagram is a sequence of type ($\dagger$) for $M$. The proof is thus complete in view of \textsc{Step~2}.
\eproof

\brem \label{rem:fg.Gtilting=>CMfinite}
By the dual of~\cite[Corollary~3.14]{ASII} about ``relative tilting modules'', a virtually Gorenstein artin algebra $\Lambda$ 
admits a finitely generated Gorenstein tilting module if and only if $\Lambda$ is of finite CM-type. Thus, the assumption on $\Lambda$ in Lemma~\ref{lem:approx.fg.tilting} is minimal. This also motivates the following problem whose answer is presently unknown to us:
\textit{Does existence of an infinitely generated Gorenstein tilting module over a virtually Gorenstein artin algebra $\Lambda$ imply that $\Lambda$ is CM-finite?}
\erem

Theorem~\ref{thm:Gtilting.class} in conjunction with Lemma~\ref{lem:approx.fg.tilting} yields the following corollary, which parallels~\cite[Proposition 4.1]{trlifaj.fpd}, and will play a key role later in Section~\ref{sec:applications}.

\bcoro \label{coro:G.tilting.class}
Let $\Lambda$ be a virtually Gorenstein ring such that $\GP_0 = \Add_\Lambda (H)$ for some $H \in \mod (\Lambda)$,
and $\mathcal{S}$ be a $\mathcal{GP}_0$-syzygy closed subclass of $\mod (\Lambda)$.
If $(\mathcal{U} , \mathcal{V})$ is the $G$-cotorsion pair cogenerated by $\mathcal{S}$, then the following statements are equivalent:
\begin{enumerate} 
\item $\mathcal{U} \subseteq \mathcal{GP}_n$ for some integer $n \geq 0$;

\item There exists a Gorenstein tilting $\Lambda$-module $T$ such that $\mathcal{V} = T^{\Gperp_\infty}$.
\end{enumerate}
If, furthermore, $\Lambda$ is an artin algebra, then $T$ is can be taken finitely generated if and only if $\mathcal{U}^\fin$ is contravariantly finite.
\ecoro

\bproof
The implication (i)$\implies$(ii) follows from Theorem~\ref{thm:Gtilting.class}: Since $\mathcal{S}$ has a representative set of elements, 
the class $\mathcal{V}$ is $G$-special preenveloping by Theorem~\ref{thm:eklof.trlifaj*}. Furthermore, every element of $\mathcal{S}$ has a 
degreewise finitely generated left $\GP_0$-resolution by Lemma~\ref{lem:H}, and this in conjunction with the fact that $\mathcal{S}$ is $\GP_0$-syzygy closed in
$\mod (\Lambda)$ implies that the class $\mathcal{V}$ is $G$-coresolving (cf. Observation~\ref{obs:*inf.classes}) as well as being 
closed under direct sums (cf. Observation~\ref{obs:BExt}-(v)) and direct summands.
It now follows from Theorem~\ref{thm:Gtilting.class}
that $\mathcal{V} = T^{\Gperp_\infty}$ for some Gorenstein tilting $\Lambda$-module $T$.
The implication (ii)$\implies$(i) follows from Lemma~\ref{lem:tilting.class}.

As for the second part of the assertion, regarding $T$ being finitely generated, it follows immediately from Lemma~\ref{lem:approx.fg.tilting}
that if $T$ is finitely generated, then $\mathcal{U}^\fin$ is contravariantly finite. Conversely, assume that $\mathcal{U}^\fin$ is contravariantly finite.
Then  for every finitely generated $\Lambda$-module $M$ there exists, by~\cite[Proposition~4.2 and Lemma~4.8]{ASI} and Remark~\ref{rem:GP.over.CM}, a $G$-exact sequence
\[\xymatrix{ o \ar[r] & M \ar[r] & V \ar[r] & U \ar[r] & o }\]
where $V \in( \mathcal{U}^\fin)^{\Gperp} \cap \mod (\Lambda)=\mathcal{V}^\fin$  and $U \in \mathcal{U}^\fin$. Thus, using an argument similar to \textsc{Step~1} and \textsc{Step~2} in the proof of Theorem~\ref{thm:Gtilting.class} we can construct a $G$-exact sequence 
\[\xymatrix{
\tau_\bullet := o \ar[r] & H \ar[r] & T_0 \ar[r] & \cdots \ar[r] & T_n \ar[r] & o
}\]
such that $T_i \in \mathcal{V}^\fin$ for all $1 \leq i \leq n$ and the cokernels of the maps in the sequence belong to $\mathcal{U}^\fin$. 
Now by a similar argument as in \textsc{Steps 3--5} in the proof of Theorem~\ref{thm:Gtilting.class}, the finitely generated $\Lambda$-module 
$T = \coprod_{i=0}^n T_i$ is Gorenstein $n$-tilting and $\mathcal{V} = T^{\Gperp_\infty}$.
\eproof

\section{Applications}
\label{sec:applications}
\label{sec:Gtilt.partial}

As a first application of Theorem~\ref{thm:Gtilting.class}, we show in the following theorem 
that every partial Gorenstein $n$-tilting module can be completed to a Gorenstein $n$-tilting module.
The problem of existence of complements to partial (standard) tilting modules was first considered by Bongartz~\cite{bongartz.partial}
in his proof of a theorem relating distinct indecomposable direct summands of a finitely generated $1$-tilting module
to the rank of the Grothendieck group of the algebra. 
It was proved \textit{loc. cit.} that every finitely generated partial $1$-tilting module over an artin algebra can be completed to a finitely generated $1$-tilting module.
Later on, and in an attempt to generalizing Bongartz's work to arbitrary finitely generated tilting modules, Rickard and Schofield~\cite{rickard.partial} proved that
a finitely generated partial $n$-tilting module over an artin algebra \emph{cannot} in general be completed to a finitely generated $n$-tilting module for $n \geq 2$.
The advent of infinitely generated tilting theory shed more light on the existence problem of complements to partial tilting modules. 
Indeed, it was proved by Angeleri-H\"{u}gel and Coelho~\cite[Theorem~2.1]{hugel.coel.partial} that 
every partial $n$-tilting module---finitely generated or not---can be completed to an $n$-tilting module. Thus, complements to partial tilting modules always exist 
in the realm of infinitely generated modules.

The relative analogue of Bongartz's result, on existence of complements to finitely generated partial $1$-tilting modules,
was proved by Auslander and Solberg in~\cite[Proposition~3.25]{ASII}. The following result should be regarded as the non-finitely generated generalization of  Auslander-Solberg result in the Gorenstein setting as well as the ``Gorenstein'' analogue of~\cite[Theorem~2.1]{hugel.coel.partial} for partial Gorenstein tilting modules.

\bthm \label{thm:partial.Gtilting}
Let $\Lambda$ be a virtually Gorenstein ring with $\GP_0 = \Add_\Lambda (H)$,  and $n \geq 0$ be an integer. 
The following statements are equivalent for a $\Lambda$-module $M$:
\begin{enumerate}
\item $M$ is a partial Gorenstein $n$-tilting module;

\item $M$ is a direct summand of a Gorenstein $n$-tilting module $T$ such that $M^{\Gperp_\infty} = T^{\Gperp_\infty}$;

\item There exists $N \in {}^{\Gperp} (M^{\Gperp_\infty})$ such that $T:=M \amalg N$ is a Gorenstein $n$-tilting module.
\end{enumerate}
\ethm

\bproof
We prove the theorem by establishing (i)$\iff$(ii) and (ii)$\iff$(iii).

Proof of (i)$\iff$(ii): Assume first that $M$ is a partial Gorenstein $n$-tilting module. 
Then the class $M^{\Gperp_\infty}$ is Gorenstein $n$-tilting by~Theorem~\ref{thm:Gtilting.class}
and so there exists a Gorenstein $n$-tilting $\Lambda$-module $U$ such that $M^{\Gperp_\infty} = U^{\Gperp_\infty}$. Since $M$ is 
partial Gorenstein $n$-tilting, $M \in {}^{\Gperp} (U^{\Gperp_\infty}) \cap U^{\Gperp_\infty} = \Add_\Lambda (U)$ by Lemma~\ref{lem:AddT}.
Therefore, there exists a cardinal $\kappa$ and a $\Lambda$-module $N$ such that $T: = U^{(\kappa)} \cong M \amalg N$.
The $\Lambda$-module $T$ is the desired Gorenstein $n$-tilting module.

Conversely, if (ii) holds, then it is readily seen that the $\Lambda$-module $M$ satisfies (G.T1) and (G.T2); see~\ref{df:G.tilting}.
Furthermore, since $M^{\Gperp_\infty}=T^{\Gperp_\infty}$, $M^{\Gperp_\infty}$ is closed under direct sums by Proposition~\ref{prop:Gtilting.class}.
Thus, $M$ is a partial Gorenstein $n$-tilting module.  

Proof of (ii)$\iff$(iii): Assume first that (ii) holds. Since $M$ is a direct summand of $T$, there exists a $\Lambda$-module $N$ such 
that $T= M \amalg N$. Then the equality $T^{\Gperp_\infty} = M^{\Gperp_\infty}$ implies $M^{\Gperp_\infty} \subseteq N^{\Gperp_\infty}$,
and so $N \in {}^{\Gperp} (N^{\Gperp_\infty}) \subseteq {}^{\Gperp} (M^{\Gperp_\infty})$.
Conversely, assume that there exists $N \in {}^{\Gperp} (M^{\Gperp_\infty})$ such that $T:=M \amalg N$ is a 
Gorenstein $n$-tilting module. Since $T = M \amalg N$, we have $T^{\Gperp_\infty} = M^{\Gperp_\infty} \cap N^{\Gperp_\infty}$.
But $N \in {}^{\Gperp} (M^{\Gperp_\infty})$ implies that $M^{\Gperp_\infty} \subseteq N^{\Gperp_\infty}$.
Therefore, $T^{\Gperp_\infty} = M^{\Gperp_\infty}$ and this complete the proof.
\eproof

We have not so far presented any example of an infinitely generated Gorenstein tilting module. In the sequel, 
we discuss results which not only provide examples of infinitely generated Gorenstein tilting modules, but also they indicate some
connections between infinitely generated Gorenstein tilting modules over virtually Gorenstein CM-finite algebras and finitistic dimension conjectures, parallel to~\cite{trlifaj.fpd}; cf. also~\cite[Chapter 17]{GT}. 

\bnotation
For a ring $\Lambda$ let $\mathcal{P}_\infty$ denote the class of all $\Lambda$-modules of finite projective dimension and $\GP_\infty$ denote the class
of all $\Lambda$-modules of finite Gorenstein projective dimension.
\enotation

\bdes{Finitistic Dimension Conjectures} \label{fpd}
Recall that given a ring $\Lambda$, 
\[ \fpd (\Lambda) :=\sup \big \{ \pd_\Lambda (M) \mid \text{$M \in \mathcal{P}_\infty$ is finitely generated} \big \} \]
is called the \emph{little finitistic dimension} of $\Lambda$, and
\[ \FPD (\Lambda) :=\sup \big \{ \pd_\Lambda (M) \mid M \in \mathcal{P}_\infty \big \} \]
is called the \emph{big finitistic dimension} of $\Lambda$.
An alternative way of computing finitistic dimensions is to use $\Gp_\infty$ and $\GP_\infty$, rather than the obvious classes $\mathcal{P}_\infty^\fin$ 
and $\mathcal{P}_\infty$. Indeed, by a well-known result of Holm~\cite[Theorem~2.28]{holm.gorenstein}
for \textit{any} ring $\Lambda$,
\begin{equation} \label{holm.FPD}
\FPD (\Lambda)  = \sup \big \{ \Gpd_\Lambda (M) \mid M \in \GP_\infty \big \} 
\end{equation}
and if $\Lambda$ is left-notherian, then we also have
\begin{equation} \label{holm.fpd}
\fpd (\Lambda)  = \sup \big \{ \Gpd_\Lambda (M) \mid \text{$M \in \GP_\infty$ is finitely generated} \big \} \: .
\end{equation}
The above formulas are sometimes more convenient to use. For example, the fact that finitistic dimension conjectures hold for Iwanaga-Gorenstein rings follows almost immediately from~\eqref{holm.FPD} and~\eqref{holm.fpd}; cf.~\cite[Corollary~12.3.2]{enochs.book}. 

Understanding the range of homological dimensions is a central problem is homological theory of rings and modules, and finitistic dimensions are natural invariants to consider in this regard. The most central problems about finitistic dimensions, known as \textit{finitistic dimension conjectures}, ask whether for an artin algebra $\Lambda$ the following statements hold:
\begin{description}
\item[(FDC1)] $\FPD (\Lambda) = \fpd (\Lambda)$;

\item[(FDC2)] $\fpd (\Lambda) < + \infty$.
\end{description}
These conjectures are due to Rosenberg and Zelinsky and were advertised by Bass~\cite[page~487]{bass.finitistic} in the 1960s. 
It was proved by Huisgen-Zimmermann~\cite{Z.domino} and Smal{\o}~\cite{smalo.fin} that (FDC1) \textit{fails} in general and so the right question to ask is for which 
algebras $\Lambda$ does the equality $\FPD (\Lambda) = \fpd (\Lambda)$ hold? The other problem, namely (FDC2), is still open in general although it is verified for
many classes of algebras. For more information about the finitistic dimension conjectures we refer to~\cite{Z.tale}. 

The connection between finitistic dimension conjectures and tilting theory was first established by Angeleri-H\"{u}gel and Trlifaj~\cite{trlifaj.fpd},
where the authors prove, among other things, that (FDC2) holds if and only if the cotorsion pair cogenerated by $\mathcal{P}_\infty^\fin$ is 
induced by a tilting module $T$ (see~\cite[Theorem 2.6]{trlifaj.fpd}), and also finitely generatedness of $T$ characterizes contravariantly finiteness of $\mathcal{P}_\infty^\fin$ over artin algebras---this is a sufficient condition proposed by Auslander and Reiten~\cite{AR.applications} for validity of finitistic dimension conjectures.
As we shall see in the sequel, there is a similar connection between Gorenstein tilting modules and finitistic dimension conjecture over
virtually Gorenstein CM-finite artin algebras.
\edes

\bnotation
When $\Lambda$ is a virtually Gorenstein ring, $(\mathcal{A}_n , \mathcal{B}_n)$ denotes the $G$-cotorsion pair cogenerated by
$\Gp_n$, and $(\mathcal{A} , \mathcal{B})$ denotes be the $G$-cotorsion pair cogenerated by $\Gp_\infty$.
\enotation

\brem
Since for every integer $n \geq 0$, $\Gp_{n} \subseteq \Gp_{n+1} \subseteq \Gp_{\infty}$,
we have the inclusions $\mathcal{A}_n \subseteq \mathcal{A}_{n+1} \subseteq \mathcal{A}$; cf.~\ref{obs:*perp}.
\erem

The first connection between finitistic dimension conjectures and Gorenstein tilting modules is indicated in the following theorem, which
states that over a virtually Gorenstein CM-finite artin algebra, the second finitistic dimension conjecture amounts to $\mathcal{B}$ being a right Gorenstein tilting class. 
It should be regarded as the ``Gorenstein analogue'' of~\cite[Theorem~2.6]{trlifaj.fpd}.

\bthm \label{thm:fpd.tilting}
If $\Lambda$ is a left noetherian virtually Gorenstein ring with $\GP_0 = \Add_\Lambda (H)$ for some $H \in \Gp_0$, 
then $\fpd_\Lambda (\Lambda) < + \infty$ if and only if  $\mathcal{B} = T^{\Gperp_\infty}$ for some Gorenstein tilting module $T$. In this case, $\fpd (\Lambda) = \Gpd_\Lambda (T)$ and $T$ can be taken $G$-properly $\Gp_\infty$-filtered.
\ethm

Before we prove the theorem, we need to prove some preparatory results regarding the basic properties of the $G$-cotorsion pairs
$(\mathcal{A}_n , \mathcal{B}_n)$ and $(\mathcal{A} , \mathcal{B})$.

\blem 
\label{Aus.Gversion}
\label{lem:(A,B)}
If $\Lambda$ is virtually Gorenstein ring, then the $G$-cotorsion pairs $(\mathcal{A} , \mathcal{B})$ and $(\mathcal{A}_n , \mathcal{B}_n)$ are complete
and $\mathcal{A}_n \subseteq \mathcal{GP}_n$.
\elem

\bproof
Since $\Gp_{\infty}$ and $\Gp_n$ (for every integer $n \geq 0$) have a representative set of elements, it follows from Theorem~\ref{thm:eklof.trlifaj*} that the $G$-cotorsion pairs $(\mathcal{A} , \mathcal{B})$ and $(\mathcal{A}_n , \mathcal{B}_n)$ are complete.  
As for the proof of the inclusion $\mathcal{A}_n \subseteq \mathcal{GP}_n$, note that by part~(ii) of Theorem~\ref{thm:eklof.trlifaj*}, 
every $M \in \mathcal{A}_n$ is a direct summand of some module $N$ which is $G$-properly $\GP_n$-filtered. 
Now since $\mathcal{GP}_n$ is closed under filtration by~\cite[Theorem~3.4]{enochs.transext} and direct summands,
$N$ and thereby $M$ belongs to $\mathcal{GP}_n$. Consequently, $\mathcal{A}_n \subseteq \GP_n$.
\eproof

\bprop \label{prop:fpd.tilting}
If $\Lambda$ is a virtually Gorenstein ring with $\GP_0 = \Add_\Lambda (H)$ for some $H \in \Gp_0$, then for every integer $n \geq 0$ the $G$-cotorsion pair $(\mathcal{A}_n , \mathcal{B}_n)$ is induced by a Gorenstein $n$-tilting $\Lambda$-module $T$. If, furthermore, $\Lambda$ is left noetherian and $n \leq \fpd (\Lambda)$, then $\Gpd_\Lambda (T)=n$.
\eprop

\bproof
Note that  the $G$-cotorsion pair $(\mathcal{A}_n , \mathcal{B}_n)$ is cogenerated by the $\GP_0$-syzygy closed subclass $\Gp_n$
of $\mod (\Lambda)$, and  $\mathcal{A}_n \subseteq \mathcal{GP}_n$ by Lemma~\ref{Aus.Gversion}. It now follows from Corollary~\ref{coro:G.tilting.class}
that $\mathcal{B}_n = T^{\Gperp_\infty}$ for some Gorenstein $n$-tilting $\Lambda$-module $T$; i.e. the $G$-cotorsion 
pair $(\mathcal{A}_n , \mathcal{B}_n)$ is induced by $T$.
Assume now that $\Lambda$ is left noetherian, $n \leq \fpd (\Lambda)$ and let $m:=\Gpd_\Lambda (T)$. By~\eqref{holm.fpd} there exists a finitely generated $\Lambda$-module $M$ with $\Gpd_\Lambda (M) = n$. Therefore, $M \in \Gp_n \subseteq \mathcal{A}_n \subseteq \mathcal{GP}_m$ by Lemma~\ref{lem:tilting.class}. Consequently, $\Gpd_\Lambda (T) =n$.
\eproof

We are now in a position to prove Theorem~\ref{thm:fpd.tilting}, which relates Gorenstein tilting modules and finitistic dimension conjectures; cf.~\cite[Theorem~2.6]{trlifaj.fpd}.

\bproof[\textsc{Proof of Theorem~\ref{thm:fpd.tilting}}]
If $n:=\fpd_\Lambda (\Lambda) < + \infty$, then $\Gp_\infty = \Gp_n$ and so the $B$-cotorsion pair
$(\mathcal{A} , \mathcal{B})=(\mathcal{A}_n , \mathcal{B}_n)$ is induced by a Gorenstein tilting module $T$ by Proposition~\ref{prop:fpd.tilting}.
Conversely, if the $G$-cotorsion pair $(\mathcal{A} , \mathcal{B})$ is induced by a Gorenstein $n$-tilting module $T$, then
$\Gp_\infty \subseteq \mathcal{A} \subseteq \mathcal{GP}_n$ by Lemma~\ref{lem:tilting.class}. That is,
$\fpd (\Lambda) \leq n < + \infty$. In this case the equality $\fpd (\Lambda) = \Gpd_\Lambda (T)$ follows from Proposition~\ref{prop:fpd.tilting}.
Furthermore, by Theorem~\ref{thm:eklof.trlifaj*}, Remark~\ref{ET.rem} and Lemma~\ref{lem:AddT}-(i),
there exists $T' \in \mathcal{A} \cap \mathcal{B} = \Add_\Lambda (T)$ such that $M:=T \amalg T'$ is $G$-properly $\Gp_\infty$-filtered. 
It is then easily seen that $M$ is a  $G$-tilting module with $\mathcal{B}=M^{\Gperp}$. 
\eproof


While Theorem~\ref{thm:fpd.tilting} relates Gorenstein tilting modules to finiteness of the little finitistic dimension, our next theorem, namely Theorem~\ref{fpd=FPD},
relates Gorenstein tilting modules to finiteness of the big finitistic dimension and it provides a sufficient condition for equality of the little and big finitistic dimensions.
We need to recall some facts about $\Sigma$-pure-injective modules before we state our next theorem.

\bdf \label{df:sigma.pure}
\label{df:pure-split}
Let $\Lambda$ be a ring. A $\Lambda$-module $T$ is called:
\begin{enumerate}
\item \emph{$\Sigma$-pure-injective} if for every cardinal $\kappa$, the $\Lambda$-module $T^{(\kappa)}$ is pure-injective;

\item \emph{$\Sigma$-pure-split} if any pure-embedding $\xymatrix{N \ar@{^{(}->}[r] & M }$ with $M \in \Add_\Lambda (T)$ splits.
\end{enumerate}
\edf

\bdes{Facts} \label{fact:sigma.pure}
The following facts about $\Sigma$-pure-injective modules will be used in the sequel:
\begin{enumerate}
\item Every finitely generated module over an artin algebra is endo-finite and hence $\Sigma$-pure-injective; see e.g.~\cite[page~4.4]{GT} 
or~\cite[Lemma~4.3]{krause.saorin}. \label{endofinite.is.spure}

\item Pure submodules of $\Sigma$-pure-injective modules are precisely their direct summands; see e.g.~\cite[Corollary 1.42]{facchini}.
In particular, every direct summand of a $\Sigma$-pure-injective module is $\Sigma$-pure-injective.

\item Every $\Sigma$-pure-injective module has a Krull-Schmidt-Azumaya decomposition; \sloppy see e.g.~\cite[Theorem~2.29]{facchini}. \label{spure.KSZ}

\item It is well-known (see e.g.~\cite[Lemma~2.32]{GT}) that every $\Sigma$-pure-injective module is $\Sigma$-pure-split. \label{spure=>split.pure}
\end{enumerate}
\edes

\bthm \label{fpd=FPD}
Let $\Lambda$ be a virtually Gorenstein ring with $\GP_0 = \Add_\Lambda (H)$ and consider the following statements about a Gorenstein $n$-tilting $\Lambda$-module $T$:
\begin{enumerate}
\item[(1)] $\Add_\Lambda (T)$ is closed under cokernels of $G$-monomorphisms;

\item[(2)] $\GP_{\infty} = {}^{\Gperp} (T^{\Gperp_\infty})$;

\item[(3)] $\mathcal{A} = {}^{\Gperp} (T^{\Gperp_\infty})$.
\end{enumerate}
Then:
\begin{enumerate}
\item The statements (1) and (2) are equivalent and they imply $\FPD (\Lambda) < + \infty$.

\item If $T$ is strongly finitely presented, then (2) implies  (3). In this case $\fpd (\Lambda) = \FPD (\Lambda) < + \infty$.

\item If $T$ is strongly finitely presented and $\Sigma$-pure-injective, then (1)--(3) are equivalent. 
\end{enumerate}
\ethm

\bproof
Part~(i): Assume (1) and notice that the inclusion ${}^{\Gperp} (T^{\Gperp_\infty}) \subseteq \GP_{\infty}$ always holds. 
As for the reverse inclusion, let $M \in \GP_{\infty}$ and take a $G$-special $T^{\Gperp_\infty}$-preenvelope
\begin{equation} \tag{$\ast$}
\xymatrix{o \ar[r] & M \ar[r] & B \ar[r] & A \ar[r] & o}
\end{equation}
of $M$. Since $M,A \in \GP_\infty$, we have $B \in T^{\Gperp_\infty} \cap \GP_\infty$.
Consequently, $B$ admits a $G$-exact $\Add_\Lambda (T)$-resolution of finite length by~\ref{prop:Gtilting.class}-(iii),
and so $B \in \Add_\Lambda (T)$ by the hypothesis. 
It then follows that  $M \in {}^{\Gperp} (T^{\Gperp_\infty})$ because $A,B \in {}^{\Gperp} (T^{\Gperp_\infty})$ by Lemma~\ref{lem:AddT}-(i) 
and the class ${}^{\Gperp} (T^{\Gperp_\infty})$ is $G$-resolving. This complete the proof of (1)$\implies$(2). As for the proof of (2)$\implies$(1),
note that $\GP_{\infty} = {}^{\Gperp} (T^{\Gperp_\infty})$ implies $\GP_{\infty} \cap T^{\Gperp_\infty} = \Add_\Lambda (T)$ by Lemma~\ref{lem:AddT}-(i). 
Since $\GP_{\infty}$ and $T^{\Gperp_\infty}$ are closed under cokernels of $G$-monomorphisms, it follows that $\Add_\Lambda (T)$ is also 
closed under cokernel of $G$-monomorphisms. Using the inclusion ${}^{\Gperp} (T^{\Gperp_\infty}) \subseteq \GP_n$, it follows from (2)
and~\eqref{holm.FPD}  that $\FPD (\Lambda) < + \infty$.

Part~(ii): Let $T$ be strongly finitely presented and assume (2). Thus, $\FPD (\Lambda) < + \infty$ which 
implies $\mathcal{A} \subseteq \GP_\infty = {}^{\Gperp} (T^{\Gperp_\infty})$; cf. Lemma~\ref{Aus.Gversion}. On the other hand, $T \in \Gp_n$ implies 
${}^{\Gperp} (T^{\Gperp_\infty}) \subseteq \mathcal{A}$. Consequently $\mathcal{A} = {}^{\Gperp} (T^{\Gperp_\infty})$.
As for the equality $\fpd (\Lambda) = \FPD (\Lambda)$, note that $\GP_\infty = {}^{\Gperp} (T^{\Gperp_\infty})$ implies $\FPD (\Lambda) \leq n$ by~\ref{lem:tilting.class} 
and $T \in \Gp_n$ implies $n \leq \fpd (\Lambda)$. Thus, $\fpd (\Lambda)=\FPD (\Lambda)$.

Part~(iii): Assume that $T$ is strongly finitely presented and $\Sigma$-pure-injective. In view of part~(ii) the proof of the assertion is complete once we
show (3)$\implies$(1). Assume (3) and notice that this implies every $G$-monomorphism in $\add_\Lambda (T)$ splits:
Indeed, if 
\[\xymatrix{\delta = o \ar[r] & A \ar[r] & B \ar[r] & C \ar[r] & o}\]
is a $G$-exact sequence with $A,B \in \add_\Lambda (T)$, then
$C \in \Gp_\infty \subseteq \mathcal{A}$ because $\add_\Lambda (T) \subseteq \Gp_\infty$ and the class $\Gp_\infty$ is closed under cokernels of monomorphisms.
On the other hand, it follows from hypothesis~(3) and Lemma~\ref{lem:AddT}-(i) that
$\add_\Lambda (T) \subseteq \mathcal{A} \cap \mathcal{B}$ and so $A \in \mathcal{B}$. Consequently, $\Gext^1_\Lambda (C,A)=o$ and 
thereby the sequence $\delta$ splits. Now we finish the proof by showing that every $G$-monomorphism in $\Add_\Lambda (T)$ is pure and hence
splits because of the fact that $T$ is $\Sigma$-pure-injective; cf.~\ref{fact:sigma.pure}. 
Let $f : A \arrow B$ be a $G$-monomorphism in $\Add_\Lambda (T)$. 
Assume without loss of generality that $f$ is an inclusion and there are elements $a_1 , \ldots , a_m \in A$ and $b_1 , \ldots , b_n \in B$ such that
for all $1 \leq i \leq m$ the equality $a_i = \sum_{j=1}^n \lambda_{ij} b_j$ holds in $B$ for some $\lambda_{ij} \in \Lambda$. 
Since $T$ is $\Sigma$-pure-injective and $A,B \in \Add_\Lambda (T)$, it follows from part~(ii) of Fact~\ref{fact:sigma.pure} that $A$ and $B$ 
are $\Sigma$-pure-injective and so they have a Krull-Schmidt-Azumaya decomposition by part~(iii) of~\ref{fact:sigma.pure}. 
Let $A = \bigoplus_{i \in I} A_i$ and $B = \bigoplus_{j \in J} B_j$ be the
Krull-Schmidt decompositions of $A$ and $B$, and note that $A_i$ and $B_i$ are finitely generated since they are direct summands of the finitely 
generated module $T$. In particular, every $A_i$ is contained in a finite direct sum of $B_j$'s. Therefore, there are finite subsets $I' \subseteq I$ and $J' \subseteq J$ such that $a_1 , \ldots , a_m \in A' := \bigoplus_{i \in I'} A_i$, $b_1 , \ldots , b_n \in B' := \bigoplus_{j \in J'} B_j$, and $A' \subseteq B'$.
The inclusion $f|_{A'} : A' \arrow B'$ in $\add_\Lambda (T)$ is a $G$-monomorphism, being a direct summand of the $G$-monomorphism $f$, and hence it
splits by the hypothesis. Consequently, there are $a'_1 , \ldots , a'_n \in A' \subseteq A$ such that
$a_i = \sum_{j=1}^n r_{ij} a'_j$ for every $1 \leq i \leq m$. This implies that $f : A \arrow B$ is a pure embedding and this completes the proof of (3)$\implies$(1).
\eproof

The Gorenstein tilting module $T$ in Theorem~\ref{thm:fpd.tilting} is \emph{not} finitely generated in general. Indeed, 
as the following theorem explains, finitely generated property of $T$ characterizes contravariantly finiteness of $\Gp_\infty$ under some mild conditions.
This result should be regarded as the Gorenstein analogue of~\cite[Proposition~4.1 and Theorem~4.2]{trlifaj.fpd}.

\bthm \label{thm:Gp.contra}
Let $\Lambda$ be a  CM-finite virtually Gorenstein artin algebra and $\mathcal{S}$ be a $\mathcal{GP}$-syzygy closed subclass of $\mod (\Lambda)$ contained $\Gp_\infty$.
If $(\mathcal{U} ,\mathcal{V})$ is the $G$-cotorsion pair cogenerated by $\mathcal{S}$ such that $\mathcal{U}^\fin$ is resolving
in $\mod (\Lambda)$, then the following statements are equivalent:
\begin{enumerate}
\item $\mathcal{U}^\fin$ is contravariantly finite in $\mod (\Lambda)$,

\item $\mathcal{V}=T^{\Gperp_\infty}$ for some finitely generated Gorenstein tilting module $T$.
\end{enumerate}
In particular, $\Gp_{\infty}$ is contravariantly finite if and only if $\mathcal{B} = T^{\Gperp_\infty}$ for some \emph{finitely generated} Gorenstein tilting module $T$.
In this case, $\FPD (\Lambda) = \fpd (\Lambda) < + \infty$.
\ethm

Before we prove Theorem~\ref{thm:Gp.contra} we require some preparatory results.

\bdes{Construction} \label{obs:A.finite}
Let $\Lambda$ be a CM-finite virtually Gorenstein artin algebra and note that  every $\Lambda$-module has a $G$-special $\mathcal{A}_n$-precover and a $G$-special $\mathcal{A}$-precover by Lemma~\ref{lem:(A,B)}. Parallel to~\cite[pages 352--353]{trlifaj.approx.little}, we now explain how to construct a $G$-special $\mathcal{A}$-precover for any $\Lambda$-module $M$ as the direct limit of a certain direct system of  $G$-special $\mathcal{A}_n$-precovers of $M$: 
Let $n \geq 0$ be an integer and find by Lemma~\ref{lem:(A,B)} a $G$-exact sequence
\[\xymatrix{ \delta_n :=  o \ar[r] & B_n \ar[r] & A_n \ar[r]^-{f_n} & M \ar[r] & o }\]
of $\Lambda$-modules where $f_n$ is a $G$-special $\mathcal{A}_n$-precover. 
Let $\xymatrixcolsep{0.58cm}\xymatrix{\gamma_n : B_n \ar[r] & B_{n+1}}$
be a $G$-special $\mathcal{B}_{n+1}$-preenvelope of $B_n$ and form the following pushout diagram with exact rows and columns.
\begin{equation} \label{eq:Af}
\begin{array}{c}
\xymatrixcolsep{0.75cm}
\xymatrixrowsep{0.75cm}
\xymatrix{
                             & o \ar[d] & o     \ar[d]     &              &    \\
\delta_n =  o \ar[r] & B_n \ar[r] \ar[d]_-{\gamma_n} & A_n \ar[r]^-{f_n} \ar[d]           & M \ar[r] \ar@{=}[d] & o \\
\delta_{n+1} =  o \ar[r] & B_{n+1} \ar[r] \ar[d]                       & A_{n+1} \ar[r]^-{f_{n+1}}  \ar[d]  & M \ar[r] & o \\
                             & \coker (\gamma_n) \ar@{=}[r] \ar[d] & \coker (\gamma_n) \ar[d]         &              &    \\
                             & o  & o          &              &   
}
\end{array}
\end{equation}
Note that the columns and rows of the diagram above are $G$-exact by~\ref{pullback.of.Bexact}.
Since $\gamma_n$ is a $G$-special $\mathcal{B}_{n+1}$-preenvelope, we have $\coker (\gamma_n) \in \mathcal{A}_{n+1}$,
and we have also $A_n \in \mathcal{A}_n \subseteq \mathcal{A}_{n+1}$. Therefore, $A_{n+1} \in \mathcal{A}_{n+1}$ by Observarion~\ref{obs:*perp}.
Consequently, the map $\xymatrix{f_{n+1}: A_{n+1} \ar[r] & M}$ is a $G$-special $\mathcal{A}_{n+1}$-precover of $M$.

Now starting from $n=0$ and following the above-mentioned recursive procedure of constructing $G$-speical $\mathcal{A}_n$-precovers of $M$, we obtain
a family 
\[\Delta : = \big \{  \xymatrix{\delta_n =  o \ar[r] & B_n \ar[r] & A_n \ar[r]^-{f_n} & M \ar[r] & o } \big \}_{n=0}^\infty \]
of $G$-exact sequences such that for every integer $n \geq 0$ there exists a commutative diagram of the form~\eqref{eq:Af}. 
The family $\Delta$ is indeed a direct system of $G$-exact sequences whose direct limit
\begin{equation} \label{eq:Af2}
\xymatrix{ \delta =  o \ar[r] & B \ar[r] & A \ar[r]^-{f} & M \ar[r] & o } \: .
\end{equation}
is $G$-exact by Lemma~\ref{lem:dirlim.of.Gexact}. As the following lemma shows, $\xymatrix{f : A \ar[r] & M}$ is indeed a $G$-special $\mathcal{A}$-precover of $M$.
\edes

\blem \label{lem:A.finite}
If $\Lambda$ is a  CM-finite virtually Gorenstein artin algebra, then:
\begin{enumerate}
\item  the $\Lambda$-homomorphism $\xymatrix{f : A \ar[r] & M}$ from~\eqref{eq:Af2}
is a $G$-special $\mathcal{A}$-precover of $M$.

\item $\mathcal{A}^\fin =\Gp_{\infty}$.
\end{enumerate}
\elem

\bproof
Part~(i): We need to show that $A \in \mathcal{A}$ and $B \in \mathcal{B}$. 
The class $\mathcal{A}$ is closed under $G$-proper $\mathcal{A}$-filtration by Relative Eklof Lemma~\ref{Eklof.lemma.Rel} and,
by Construction~\ref{obs:A.finite}, the $\Lambda$-module $A$ is  $G$-properly $\mathcal{A}$-filtered by the family $\{A_n\}_{n=0}^\infty$.
Therefore, $A \in \mathcal{A}$. On the other hand, $\Gp_{\infty} = \bigcup_{n=0}^\infty \Gp_n$ implies $\mathcal{B} = \bigcap_{n=0}^\infty \mathcal{B}_n$.
Therefore, $B \in \mathcal{B}$ amounts to $B \in \mathcal{B}_n$ for every integer $n \geq 0$. But, keeping the notations of~\ref{obs:A.finite} in mind,
it is easy to see that for every integer $n \geq 0$, $B = \bigcup^\infty_{m=0} B_m =\bigcup^\infty_{m \geq n} B_m \in \mathcal{B}_n$.
The proof of part~(i) is thus complete.

Part~(ii): Clearly, $\Gp_\infty \subseteq \mathcal{A}^\fin$. As for the reverse inclusion, let $M \in \mathcal{A}^\fin$ 
and note that by part~(i) of the lemma the $G$-exact sequence~\eqref{eq:Af2} associated with $M$ splits. Consequently,
$M$ is a direct summand of $A = \bigcup_{n=0}^{\infty} A_n$. Since $M$ is finitely generated, there exists some integer $n \geq 0$
such that $M$ is a direct summand of $A_n$. Thus, $M \in \mathcal{A}_n \subseteq \GP_n$ in view of Lemma~\ref{lem:(A,B)},
and therefore $M \in \Gp_{\infty}$. This completes the poof of part~(ii).
\eproof

\bproof[\textsc{Proof of Theorem~\ref{thm:Gp.contra}}]
By Corollary~\ref{coro:G.tilting.class}, it suffices to show that $\mathcal{U} \subseteq \GP_n$ for some integer $n \geq 0$.
Since $\mathcal{S} \subseteq \Gp_{\infty}$, we have $\mathcal{U} \subseteq \mathcal{A}$ and hence
$\mathcal{U}^\fin \subseteq \Gp_{\infty}$ by Lemma~\ref{lem:A.finite}. Consequently, there exists  
an integer $n \geq 0$ such that $\mathcal{U}^\fin \subseteq \Gp_n$ by~\cite[Corollary~3.10]{AR.applications}.
Now by the inclusion $\mathcal{S} \subseteq \mathcal{U}^\fin$, every $\Lambda$-module in $\mathcal{U}$ is a direct summand of a 
$G$-properly $\GP_n$-filtered module by Theorem~\ref{thm:eklof.trlifaj*} and hence it belongs to $\GP_n$ by Lemma~\ref{Aus.Gversion}.
That is, $\mathcal{U} \subseteq \GP_n$ which finishes the proof.
\eproof

\bdes{Final Remark}
Contravariantly finiteness of the class $\Gp_\infty$ in Theorem~\ref{thm:Gp.contra} and its relation to finitistic dimension conjectures has also been considered in~\cite[Proposition 4.8]{xiIII}, and it is indeed reminiscent of ``contravariantly finiteness of the class $\mathcal{P}_\infty^\fin$''---a sufficient condition for validity of
the finitistic dimension conjectures due to Auslander and Reiten~\cite{AR.applications}. It is very tempting in the first place to think that ``contravariantly finiteness of $\Gp_\infty$'' is a new sufficient condition for validity of finitistic dimension conjectures. But one should resists this temptation as the first author and Jan \v{S}aroch have proved recently in an as yet unpublished work that over artin algebras ``contravariantly finiteness of the class $\Gp_\infty$'' implies  ``contravariantly finiteness of the class $\mathcal{P}_\infty^\fin$'', and the converse holds when the class $\Gp_0$ is contravariantly finite. 
\edes


\section*{Acknowledgements}
Part of this work was done during the first author's visit to the Department of Algebra at Charles University, Czech Republic.
He would like to thank University of Tehran for financial support of the visit and express his gratitude to Jan Trlifaj for his warm hospitality.
The research of Siamak Yassemi is supported by Iran National Science Foundation (INSF), Project No. 95831492.

\bibliographystyle{amsplain} 
\bibliography{thebib}

\providecommand{\bysame}{\leavevmode\hbox to3em{\hrulefill}\thinspace}
\providecommand{\MR}{\relax\ifhmode\unskip\space\fi MR }
\providecommand{\MRhref}[2]{%
  \href{http://www.ams.org/mathscinet-getitem?mr=#1}{#2}
}
\providecommand{\href}[2]{#2}
\begin{thebibliography}{10}

\bibitem{hom.flat}
Stephen~T Aldrich, Edgar~E Enochs, and Juan~A Lopez~Ramos, \emph{Derived
  functors of {H}om relative to flat covers}, Mathematische Nachrichten
  \textbf{242} (2002), no.~1, 17--26.

\bibitem{hugel.coel}
Lidia Angeleri-H\"{u}gel and Fl{\'a}vio~Ulhoa Coelho, \emph{Infinitely
  generated tilting modules of finite projective dimension}, Forum
  Mathematicum, vol.~13, Berlin; New York: De Gruyter, c1989-, 2001,
  pp.~239--250.

\bibitem{hugel.coel.partial}
Lidia Angeleri-H{\"u}gel and Fl{\'a}vio~Ulhoa Coelho, \emph{Infinitely
  generated complements to partial tilting modules}, Mathematical Proceedings
  of the Cambridge Philosophical Society, vol. 132, Cambridge University Press,
  2002, pp.~89--96.

\bibitem{tilting.handbook}
Lidia Angeleri-H{\"u}gel, Dieter Happel, and Henning Krause, \emph{Handbook of
  tilting theory}, vol.~13, Cambridge University Press, 2007.

\bibitem{tilting.preenv}
Lidia Angeleri-H{\"u}gel, Alberto Tonolo, and Jan Trlifaj, \emph{Tilting
  preenvelopes and cotilting precovers}, Algebras and Representation Theory
  \textbf{4} (2001), no.~2, 155--170.

\bibitem{trlifaj.fpd}
Lidia Angeleri-H\"{u}gel and Jan Trlifaj, \emph{Tilting theory and the
  finitistic dimension conjectures}, Transactions of the American Mathematical
  Society \textbf{354} (2002), no.~11, 4345--4358.

\bibitem{AR.applications}
{M}aurice Auslander and Idun {R}eiten, \emph{Applications of contravariantly
  finite subcategories}, Advances in Mathematics \textbf{86} (1991), no.~1,
  111--152.

\bibitem{aus.preproj}
Maurice Auslander and Sverre~O Smal{\o}, \emph{Preprojective modules over
  {A}rtin algebras}, Journal of algebra \textbf{66} (1980), no.~1, 61--122.

\bibitem{ASI}
Maurice Auslander and {\O}yvind Solberg, \emph{Relative homology and
  representation theory {I}. {R}eative homology and homologically finite
  subcategories}, Communications in Algebra \textbf{21} (1993), no.~9,
  2995--3031.

\bibitem{ASII}
\bysame, \emph{Relative homology and representation theory {II}: {R}elative
  cotilting theory}, Communications in Algebra \textbf{21} (1993), no.~9,
  3033--3079.

\bibitem{ASIII}
\bysame, \emph{Relative homology and representation theory {III}: {C}otilting
  modules and {W}edderburn correspondence}, Communications in Algebra
  \textbf{21} (1993), no.~9, 3081--3097.

\bibitem{bass.finitistic}
Hyman Bass, \emph{Finitistic dimension and a homological generalization of
  semi-primary rings}, Transactions of the American Mathematical Society
  \textbf{95} (1960), no.~3, 466--488.

\bibitem{beligiannis.CM}
Apostolos Beligiannis, \emph{Cohen-{M}acaulay modules,(co)torsion pairs and
  virtually {G}orenstein algebras}, Journal of Algebra \textbf{288} (2005),
  no.~1, 137--211.

\bibitem{beligiannis.CMtype}
\bysame, \emph{On algebras of finite {C}ohen--{M}acaulay type}, Advances in
  Mathematics \textbf{226} (2011), no.~2, 1973--2019.

\bibitem{belig.krause.thick}
Apostolos Beligiannis, Henning Krause, et~al., \emph{Thick subcategories and
  virtually {G}orenstein algebras}, Illinois Journal of Mathematics \textbf{52}
  (2008), no.~2, 551--562.

\bibitem{BR}
Apostolos Beligiannis and Idun Reiten, \emph{Homological and homotopical
  aspects of torsion theories}, Mem. Amer. Math. Soc., 2007.

\bibitem{bongartz.partial}
Klaus Bongartz, \emph{Tilted algebras}, Representations of Algebras (Berlin,
  Heidelberg) (Maurice Auslander and Emilo Lluis, eds.), Springer Berlin
  Heidelberg, 1981, pp.~26--38.

\bibitem{brenner.butler}
Sheila Brenner and Michael~CR Butler, \emph{Generalizations of the
  {B}ernstein-{G}elfand-{P}onomarev reflection functors}, Representation theory
  II, Springer, 1980, pp.~103--169.

\bibitem{chen.balanced}
Xiao-Wu Chen, \emph{Homotopy equivalences induced by balanced pairs}, Journal
  of Algebra \textbf{324} (2010), no.~10, 2718--2731.

\bibitem{chen.gorenstein}
\bysame, \emph{Gorenstein homological algebra of {A}rtin algebras}, arXiv
  preprint arXiv:1712.04587 (2017).

\bibitem{inf.tilting}
Riccardo Colpi and Jan Trlifaj, \emph{Tilting modules and tilting torsion
  theories}, Journal of Algebra \textbf{178} (1995), no.~2, 614--634.

\bibitem{tilting.sub}
Zhenxing Di, Jiaqun Wei, Xiaoxiang Zhang, and Jianlong Chen, \emph{Tilting
  subcategories with respect to cotorsion triples in abelian categories},
  Proceedings of the Royal Society of Edinburgh Section A: Mathematics
  \textbf{147} (2017), no.~4, 703--726.

\bibitem{EM}
Samuel Eilenberg and John~C Moore, \emph{Foundations of relative homological
  algebra}, no.~55, American Mathematical Soc., 1965.

\bibitem{eklof}
Paul~C Eklof, \emph{Homological algebra and set theory}, Transactions of the
  American Mathematical Society \textbf{227} (1977), 207--225.

\bibitem{trlifaj.Extvanish}
Paul~C Eklof and Jan Trlifaj, \emph{How to make {E}xt vanish}, Bulletin of the
  London Mathematical Society \textbf{33} (2001), no.~1, 41--51.

\bibitem{enochs.coversenv}
Edgar~E Enochs, \emph{Injective and flat covers, envelopes and resolvents},
  Israel Journal of Mathematics \textbf{39} (1981), no.~3, 189--209.

\bibitem{enochs.Fbalance}
\bysame, \emph{Balance with flat objects}, Journal of Pure and Applied Algebra
  \textbf{219} (2015), no.~3, 488--493.

\bibitem{enochs.transext}
Edgar~E Enochs, Alina Iacob, Overtoun~MG Jenda, et~al., \emph{Closure under
  transfinite extensions}, Illinois Journal of Mathematics \textbf{51} (2007),
  no.~2, 561--569.

\bibitem{enochs.balanced}
Edgar~E Enochs and Overtoun~MG Jenda, \emph{Balanced functors applied to
  modules}, Journal of Algebra \textbf{92} (1985), no.~2, 303--310.

\bibitem{enochs.Gext}
\bysame, \emph{Gorenstein balance of {H}om and tensor}, Tsukuba journal of
  mathematics \textbf{19} (1995), no.~1, 1--13.

\bibitem{enochs.book}
E.E. Enochs and O.M.G. Jenda, \emph{Relative {H}omological {A}lgebra}, De
  Gruyter Expositions in Mathematics, De Gruyter, 2011.

\bibitem{facchini}
Alberto Facchini, \emph{Module {T}heory: {E}ndomorphism rings and direct sum
  decompositions in some classes of modules}, Modern Birkh{\"a}user Classics,
  Springer Basel, 2012.

\bibitem{GT}
R{\"u}diger G{\"o}bel and Jan Trlifaj, \emph{Approximations and {E}ndomorphism
  {A}lgebras of {M}odules: {V}olume 1--{A}pproximations/{V}olume
  2--{P}redictions}, vol.~41, Walter de Gruyter, 2012.

\bibitem{happel.book}
Dieter Happel, \emph{Triangulated categories in the representation of finite
  dimensional algebras}, vol. 119, Cambridge University Press, 1988.

\bibitem{tiltedalg.happel}
Dieter Happel and Claus~Michael Ringel, \emph{Tilted algebras}, Transactions of
  the American Mathematical Society (1982), 399--443.

\bibitem{holm.Gext}
Henrik Holm, \emph{Gorenstein derived functors}, Proceedings of the American
  Mathematical Society \textbf{132} (2004), no.~7, 1913--1923.

\bibitem{holm.gorenstein}
\bysame, \emph{Gorenstein homological dimensions}, Journal of pure and applied
  algebra \textbf{189} (2004), no.~1, 167--193.

\bibitem{holm.rel.Ext}
\bysame, \emph{Relative {E}xt groups, resolutions, and {S}chanuel classes},
  Osaka J. Math. \textbf{45} (2008), no.~3, 719--735.

\bibitem{Z.domino}
Birge~Zimmermann Huisgen, \emph{Homological domino effects and the first
  finitistic dimension conjecture}, Inventiones mathematicae \textbf{108}
  (1992), no.~1, 369--383.

\bibitem{Z.tale}
\bysame, \emph{The {F}initistic {D}imension {C}onjectures---{A} {T}ale of 3.5
  {D}ecades}, pp.~501--517, Springer Netherlands, Dordrecht, 1995.

\bibitem{krause.saorin}
Henning Krause and Manuel Saor{\'\i}n, \emph{On minimal approximations of
  modules}, Trends in the Representation Theory of Finite Dimensional Algebras
  (Edward~L. {G}reen and Birge {H}uisgen {Z}immermann, eds.), Contemporary
  mathematics - American Mathematical Society, vol. 229, American Mathematical
  Society, 1998, pp.~227--236.

\bibitem{rel.cotorsion}
Huan~Huan Li, JunFu Wang, and ZhaoYong Huang, \emph{Applications of balanced
  pairs}, Science China Mathematics \textbf{59} (2016), no.~5, 861--874.

\bibitem{miyashita}
Yoichi Miyashita, \emph{Tilting modules of finite projective dimension},
  Mathematische Zeitschrift \textbf{193} (1986), no.~1, 113--146.

\bibitem{rickard.partial}
Jeremy Rickard and Aidan Schofield, \emph{Cocovers and tilting modules},
  Mathematical Proceedings of the Cambridge Philosophical Society \textbf{106}
  (1989), no.~1, 1–5.

\bibitem{RotmanSE}
Joseph Rotman, \emph{An introduction to homological algebra}, Springer Science
  \& Business Media, 2008.

\bibitem{salce}
Luigi Salce, \emph{Cotorsion theories for abelian groups}, Symposia Math,
  vol.~23, 1979, p.~3.

\bibitem{smalo.fin}
Sverre~O Smal{\o}, \emph{Homological differences between finite and infinite
  dimensional representations of algebras}, Infinite length modules, Springer,
  2000, pp.~425--439.

\bibitem{trlifaj.approx.little}
Jan Trlifaj, \emph{Approximations and the little finitistic dimension of
  artinian rings}, Journal of Algebra \textbf{246} (2001), no.~1, 343--355.

\bibitem{trlifaj.infinite.cotorsion}
\bysame, \emph{Infinite dimensional tilting modules and cotorsion pairs},
  LONDON MATHEMATICAL SOCIETY LECTURE NOTE SERIES \textbf{332} (2007), 279.

\bibitem{wei}
Jiaqun Wei, \emph{A note on relative tilting modules}, Journal of Pure and
  Applied Algebra \textbf{214} (2010), no.~4, 493--500.

\bibitem{xiIII}
Changchang Xi, \emph{On the finitistic dimension conjecture, iii: Related to
  the pair eae⊆a}, Journal of Algebra \textbf{319} (2008), no.~9, 3666 --
  3688.

\bibitem{Gtilting}
Liang Yan, Weiqing Li, and Baiyu Ouyang, \emph{Gorenstein {C}otilting and
  {T}ilting {M}odules}, Communications in Algebra \textbf{44} (2016), no.~2,
  591--603.

\end{thebibliography}
\end{document}